\documentclass{jnmp}
\usepackage{amsmath}
\JNMPnumberwithin{equation}{section}
\def\R{{\mathbb R}}

\begin{document}
\renewcommand{\evenhead}{A~M Grundland and W~J~Zakrzewski}
\renewcommand{\oddhead}{On $CP^{1}$ and $CP^{2}$ Maps and Weierstrass
 Representations}

\thispagestyle{empty}

\FirstPageHead{10}{1}{2003}{\pageref{grundland-firstpage}--\pageref{grundland-lastpage}}{Article}

\copyrightnote{2003}{A~M Grundland and W~J~Zakrzewski}

\Name{On $\boldsymbol{CP^{1}}$ and $\boldsymbol{CP^{2}}$ Maps and Weierstrass
 Representations for Surfaces Immersed\\ into Multi-Dimensional Euclidean Spaces}
\label{grundland-firstpage}

\Author{A~M GRUNDLAND~$^\dag$ and W~J~ZAKRZEWSKI~$^\ddag$}

\Address{$^\dag$~Centre de Recherches Math{\'e}matiques, Universit{\'e} de Montr{\'e}al,\\
~~C.\ P.~6128, Succ.\ Centre-ville, Montr{\'e}al, (QC) H3C 3J7, Canada\\
~~E-mail: Grundlan@crm.umontreal.ca\\[10pt]
$^\ddag$~Department of Mathematical Sciences, University of Durham,
 Durham DH1 3LE, UK\\
~~E-mail: W.J.Zakrzewski@durham.ac.uk}

\Date{Received June 07, 2002; Accepted July 24, 2002}

\begin{abstract}
\noindent
An extension of the classic Enneper--Weierstrass representation
 for conformally pa\-ra\-met\-ri\-sed surfaces in multi-dimensional spaces
 is presented. This is based on low dimensional $CP^1$
 and $CP^2$ sigma models which allow the study of the constant mean
 curvature (CMC) surfaces immersed into Euclidean 3- and
 8-dimensional spaces, respectively. Relations of Weierstrass
 type systems to the equations of these sigma models
  are established.
 In  particular, it is demonstrated that the generalised Weierstrass
 representation can admit different CMC-surfaces in ${\R}^3$ which
 have globally the same Gauss map. A new procedure for constructing
 CMC-surfaces in ${\R}^n$ is presented and illustrated in some explicit
 examples.
 \end{abstract}

 \section{Introduction}

   In this paper we study two-dimensional surfaces
 conformally immersed into multi dimensional spaces with Euclidean
 metric. We present explicit formulae for the position vector $X:
 {\mathcal D} \in {\mathbb C} \rightarrow {\R}^3$ of a
 surface for which $X$ satisfies the Gauss--Weingarten and Gauss--Codacci
 equations. Such formulae describing minimal surfaces (i.e.\
 zero mean curvature $H=0$) imbedded in three-dimensional space were
 first formulated by Enneper~\cite{Enneper} and Weierstrass~\cite{Weierstrass} about one and
 half century ago. These authors  considered two holomorphic
 functions $\psi(z)$ and
 $\phi(z)$ of a complex variable $z \in {\mathbb C}$ and introduced a three
 component complex vector valued function $w=(w_1,w_2,w_3): {\mathcal D}
 \rightarrow {\mathbb C}^3$ which is required to satisfy the following
 linear differential equations
 \begin{gather}
 \partial w_1 = \frac{1}{2} \left(\psi^{2} - \phi^{2}\right),
 \qquad
 \partial w_2 = \frac{i}{2}\left(\psi^{2} + \phi^{2}\right),
 \qquad
 \partial w_3= - \psi \phi, \nonumber\\
 \bar{\partial} \psi = 0, \qquad \bar{\partial} \phi = 0,\label{1.1}
 \end{gather}
 where the derivatives are abbreviated $\partial = \partial/ \partial z$
 and the bar denotes the complex congugation. Then they showed that
 the real
 vector valued functions
 \begin{equation}\label{1.2}
 X = \left( {\rm Re} \int_{0}^{z} \frac{1}{2}\left(\psi^{2} - \phi^{2}\right)\,dz',
  {\rm Re} \int_{0}^{z} \frac{i}{2}\left(\psi^{2} + \phi^{2}\right)\,dz',
 -  {\rm Re} \int_{0}^{z} \psi  \phi\,dz'\right)
 \end{equation}
 can be considered as components of a position vector
  of a minimal surface, immersed
 into ${\R}^3$, with the conformal metric
 \begin{equation}\label{1.3}
 ds^2 = \left(|\psi|^{2} + |\phi|^{2}\right)\sp2\,dz\,d \bar{z},
 \end{equation}
 where $z$ and $\bar{z}$ are local coordinates on ${\cal D}$ and
 the coordinate lines $z={\rm const}$, and $\bar z={\rm const}$ describe
 geodesics on this surface. Since then
 this
 idea has been developed by many authors,  for a review of the subject
 see e.g.~\cite{book1,K,book2} and references therein.
   The theory of minimal, or in general, constant mean curvature
 (CMC)-surfaces plays an
 important
 role in several
 applications to problems arising both in mathematics and in physics. In
 particular, many interesting applications can be found in such
 diverse areas of physics as: the fields of two-dimensional
gravity~\cite{book3,Carroll}, string theory~\cite{Land}, quantum field
 theory~\cite{Vis,book3}, statistical physics~\cite{book5,book6}, fluid
dynamics~\cite{book7}, theory of fluid membranes~\cite{Z,book5}.
 One interesting application involves the Canham--Helfrich
 membrane model~\cite{Can,Hel}.
  This model can explain some basic features and equilibrium
 shapes both for biological membranes and liquid interfaces~\cite{book8}.

   Our approach involves modifying the original
 Enneper--Weierstrass representation (\ref{1.2}) by adding to it extra
 terms. For this purpose it is convenient to exploit the
 connection
 between Weierstrass systems, $CP^1$ and $CP^2$ sigma model equations, and
 their Lax representations. We demonstrate that, through these links,
 conformal immersion of CMC-surfaces into 3- and 8-dimensional spaces can
 be formulated. We show that a large classes of solutions of the
 Weierstrass
 system can be obtained and, consequently, can provide new classes of
 conformally parametrised CMC-surfaces in multi-dimensional spaces.

   The paper is organized as follows.
 In Section~2, we rederive the classical Enneper--Weierstrass
 representation for minimal surfaces immersed into ${\R}^3$.
 In the next section we
 describe, in detail, the generalised Weierstrass formulae for
 CMC-surfaces
 into ${\R}^3$ in the context of the $CP^1$ sigma model and discuss
 some geometric aspects of $CP^1$ maps. The following section deals with
 $CP^2$ maps and the corresponding Weierstrass representation for
 conformally
 parametrised surfaces immersed into ${\R}^8$. It also presents some
 geometric
 characteristics of CMC-surfaces.
  In Section~5 our theoretical
 considerations are illustrated by explicit examples and new interesting
 CMC-surfaces are found. The last section presents final remarks and
 mentions possible future developments.

 \section{The Enneper--Weierstrass formulae\\ for
  minimal surfaces in $\boldsymbol{{\R}^{3}}$}

 Let ${\mathbb M}^2$ be a smooth orientable surface in $3$-dimensional Euclidean
 space ${\R}^3$.
 The surface~${\mathbb M}^2$ is described by a real vector-valued function
 \begin{equation}\label{2.1}
 X = (X_{1}, X_{2}, X_{3}) : D \rightarrow  {\R}^3,
 \end{equation}
 where $D$ is a region in the complex plane ${\mathbb C}$. The metric
 is assumed to be conformally flat
 \begin{equation}\label{2.2}
 ds^2 = e^{2u} \, dz \, d \bar{z}
 \end{equation}
 for any real valued function $u$ of $z$ and $\bar{z}$. The conformal
 parametrisation of the surface ${\mathbb M}^2$ implies the following normalization
 of
 the position vector $X(z, \bar{z})$
 \begin{equation}\label{2.3}
 ( \partial X, \partial X) = 0,\qquad
 ( \partial X, \bar{\partial} X) = \frac{1}{2} e^{2u},
 \end{equation}
 where  the brackets $(\cdot ,\cdot )$
 denote the standard scalar product in ${\R}^3$. The
 tangent vectors $\partial X$ and $\bar{\partial} X$ and the real
 unit normal vector $N$ on the surface ${\mathbb M}^2$ satisfy the obvious
  relations
 \begin{equation}\label{2.4}
 (\partial X, N) = 0, \qquad
 (N, N) = 1.
 \end{equation}

 Equations of a moving complex frame $\xi = ( \partial X, \bar{\partial}
 X, N)^{T}$
 satisfy the following Gauss--Weingarten equations (see e.g.~[5])
 \begin{equation}\label{2.6}
 \partial \xi = U \xi, \qquad
 \bar{\partial} \xi = V \xi,
 \end{equation}
 where $3\times 3$ matrices $U$ and $V$ have the form
 \begin{equation}\label{2.7}
 U = \left(
 \begin{array}{ccc}
 2 \partial u  &  0  &   J   \\
 0      &   0   &  \frac{1}{2} H e^{2u}   \\
 -H      & -2 e^{-2u}  J  &   0  \\
 \end{array}   \right),   \qquad
 V=  \left(
 \begin{array}{ccc}
 0  &   0   &  \frac{1}{2}  H e^{2u}   \\
 0  & 2 \bar{\partial} u   &  \bar{J}  \\
 -2 e^{-2u}  \bar{J}  &  - H    & 0    \\
 \end{array}  \right),
 \end{equation}
 and the following notation has been introduced
 \begin{equation}\label{2.8}
 J = (\partial^2 X, N),   \qquad  H = 2 e^{-2u}
 (\partial \bar{\partial} X, N).
 \end{equation}
 Formulae (\ref{2.6}) are compatible with the scalar
 products (\ref{2.3}) and (\ref{2.4}). From (\ref{2.6}) and~(\ref{2.7}) we
 can derive
 the equation for the unit normal vector $N$
 \begin{equation}\label{2.9}
 \partial \bar{\partial} N + (\partial N,
 \bar{\partial} N ) N + \bar{\partial} H
 \partial X + \partial H \bar{\partial} X= 0.
 \end{equation}
 The corresponding Gauss--Codazzi equations of the
 conformally parametrised surface ${\mathbb M}^2\!\subset {\R}^3$
 are the compatibility conditions of equations (\ref{2.6})
 and have the following form
 \begin{gather}
\partial \bar{\partial} u + \frac{1}{4} H^2 e^{2u}
 - 2 |J|^2 e^{-2u}  = 0, \label{2.10}\\
\bar{\partial} J = \frac{1}{2} \partial H e^{2u},
 \qquad \partial \bar{J} = \frac{1}{2} \bar{\partial} H e^{2u}.
 \label{2.11}
 \end{gather}

 The aim of this section is to rederive the original
 Enneper--Weierstrass formulae \cite{Enneper,Weierstrass} for inducing
 minimal surfaces in ${\R}^3$. For surfaces with $H=0$
  the formulae given above simplify considerably. We
 focus our attention
 on the construction of the
 explicit formula for the position vector $X(z, \bar{z})$ of
 conformally parametrised surfaces into ${\R}^3$ for which
 equations (\ref{2.3}), (\ref{2.10}) and (\ref{2.11}) are fulfilled.

 For computational purposes, it is useful to examine equations
 (\ref{2.3}), (\ref{2.10}) and (\ref{2.11}) in terms of a
 two-component object which, in fact, is a spinor, but
 its spinorial nature is not relevant to our discussion:
 $\phi = (\psi_{1}, \psi_{2})^{T} \in {\mathbb C}^2$.
  We show that, by quadratures, we can determine  the
 coordinates of the position vector $X(z, \bar{z})$ in terms
 of the components of $\phi$ satisfying equations
 (\ref{2.3}) and (\ref{2.10})--(\ref{2.11}).

Let us consider the complex vector $\vec{w}$ in ${\mathbb C}^3$
 equal to one of the tangent vectors, say,  $\partial X$
 \begin{equation}\label{2.12}
 \vec{w} = ( w_{1}, w_{2}, w_{3}) = \partial X,
 \qquad w_{i} \in {\mathbb C}, \quad i=1,2,3,
 \end{equation}
 the $2$ by $2$ traceless matrix
 \begin{equation}\label{2.13}
 w = \left(
 \begin{array}{cc}
 w_{3}  &   w_{1}  - i w_{2}   \\
 w_{1} + i w_{2}  &   - w_{3}  \\
 \end{array}   \right),
 \qquad    \mbox{tr}\, w  = 0,
 \end{equation}
 and the map
 \begin{equation}\label{2.14}
  \vec{w} : {\mathbb C}  \rightarrow w = w_{i} \sigma_{i}\in sl(2,{\mathbb C}),
 \end{equation}
 where $\sigma_{i}$ are the Pauli matrices
 \begin{equation}\label{2.15}
 \sigma_{1}  = \left(
 \begin{array}{cc}
 0 & 1  \\
 1 & 0   \\
 \end{array}   \right), \qquad
 \sigma_{2}  = \left(
 \begin{array}{cc}
 0 & -i  \\
 i &  0  \\
 \end{array}   \right), \qquad \sigma_{3}  = \left(
 \begin{array}{cc}
 1 & 0   \\
 0 & -1  \\
 \end{array}    \right).
 \end{equation}
 The map (\ref{2.14}) satisfies
 \begin{equation}\label{2.16}
 \vec{w}^2 = - \det  w.
 \end{equation}
 From (\ref{2.16}), the determinant of the matrix $w$ vanishes
 if and only if the vector $\vec{w}$ is null, which
 coincides with the first condition in (\ref{2.3}).
 Hence, using (\ref{2.13}), we can express uniquely
 the null vector $\vec{w}$ in terms of the complex
 two-component vector $\phi$ as follows
 \begin{equation}\label{2.17}
 w_{1} = \frac{1}{2} \left( \psi_{1}^{2} - \psi_{2}^{2}\right),
 \qquad
 w_{2} = \frac{i}{2} \left( \psi_{1}^{2} + \psi_{2}^{2}\right),
 \qquad
 w_{3} = - \psi_{1} \psi_{2}.
 \end{equation}
 From the assumption (\ref{2.12}) that the null vector $\vec{w}$ is
 equal to the tangent vector $\partial X$, we can express
 $\partial X$ in terms of
 $\psi_{1}$ and $\psi_{2}$ as follows:
 \begin{equation}\label{2.18}
 \partial X_{1} = \frac{1}{2} \left( \psi_{1}^2 - \psi_{2}^2\right),
 \qquad
 \partial X_{2} = \frac{i}{2} ( \psi_{1}^2 + \psi_{2}^2),
 \qquad
 \partial X_{3} = - \psi_{1} \psi_{2},
 \end{equation}
 which coincide with expression (\ref{1.1}).
 The Enneper--Weierstrass representation for surfaces in ${\R}^3$ are
 obtained under the additional assumption that  $\psi_{1}$
 and $\psi_{2}$ are arbitrary holomorphic functions
 of the complex variable $ z \in {\mathbb C}$.
 Then, integrating equations (\ref{2.18}) and taking into account the
 reality
 condition of the position vector
 \begin{equation}\label{2.19}
 X(z, \bar{z}) = \bar{X} ( z, \bar{z}),
 \end{equation}
 we obtain the following representation \cite{Weierstrass}
 \begin{gather}
X_{1} = \frac{1}{2} \int_{0}^{z}
 \left( \psi_{1}^{2} - \psi_{2}^2\right)  dz'
 + \frac{1}{2} \int_{0}^{\bar{z}}
\left( \bar{\psi}_{1}^2 - \bar{\psi}_{2}^2\right)  dz',
 \nonumber\\
X_{2} = \frac{i}{2} \int_{0}^{z}
 \left(\psi_{1}^2 + \psi_{2}^{2}\right)  dz'
 - \frac{i}{2} \int_{0}^{\bar{z}}
 \left( \bar{\psi}_{1}^2 + \bar{\psi}_{2}^2\right)  d \bar{z}',
 \nonumber\\
X_{3} = - \int_{0}^{z} \psi_{1} \psi_{2} \,dz'
 - \int_{0}^{\bar{z}} \bar{\psi}_{1} \bar{\psi}_{2} \, d \bar{z},
 \label{2.20}
 \end{gather}
 which, in fact, is equivalent to (1.2).
 Next, from (\ref{2.20}) and invoking the second condition~(\ref{2.3}) we find that
 \begin{equation}\label{2.21}
 u = \ln \left(|\psi_{1}|^2 + |\psi_{2}|^2\right).
 \end{equation}
 Substituting (\ref{2.21}) into the Gauss--Codazzi equations
 (\ref{2.10})--(\ref{2.11}),
 we obtain
 \begin{equation}\label{2.23}
 \partial \bar{\partial} \ln p^2
 -2 |J|^2 p^{-2} = 0,
 \end{equation}
 where
 \begin{equation}\label{2.22}
 p= |\psi_{1}|^2 + |\psi_{2}|^2.
 \end{equation}
  By virtue of (\ref{2.20}),
 we find that $J$ and $H$ defined by (2.7),
 when expressed in terms of $\psi_1$ and $\psi_2$ become
 \begin{equation}\label{2.24}
 J = \psi_{1} \partial \psi_{2} - \psi_{2} \partial \psi_{1},\qquad H=0
 \end{equation}
 and $J$ is analytic, i.e.\  $
 \bar{\partial} J = 0$.

  Note that the
 direction of $\phi = ( \psi_{1}, \psi_{2})^{T}$
 is arbitrary, but its length is fixed by (\ref{2.23}).
 Note also that after the change of variable
 $\varphi = 2 \ln p$
 equation (\ref{2.23}) becomes
 \begin{equation}\label{2.26}
 \partial \bar{\partial} \varphi
 = 2  |J|^{2} e^{- \varphi},
 \qquad
 \bar{\partial} J = 0.
 \end{equation}

 \section{The generalised Weierstrass formulae\\
 for CMC-surfaces in  $\boldsymbol{{\R}^3}$}

 The Wierstrass--Enneper formulae for inducing minimal surfaces, and their
 generalisations,
 has been studied for a long time by many authors
 (e.g.\ \cite{book1,book9,book12} and references therein). This topic
 has most recently been treated by B~Konopelchenko and I~Taimanov~\cite{Tai}.
 In this paper, Konopelchenko and  Taimanov, established a direct
 connection between certain classes of
 CMC-surfaces and an integrable finite-dimensional Hamiltonian system.
 For a summary of their results, see~\cite{Land2}. There it is shown
 that
 to any  solution $\phi = (\psi_{1}, \psi_{2})^{T}$ of the first
 order equations (which resemble Dirac equations)
 \begin{equation}\label{3.1}
 \partial \psi_{1} = p \psi_{2},\qquad  \bar{\partial} \psi_{2} = - p
 \psi_{1}, \qquad  p  = |\psi_{1}|^2 + | \psi_{2}|^2,
 \end{equation}
 one can associate a CMC-surface immersed into ${\R}^3$ with radius
 vector $X(z, \bar{z})$ of the form~(\ref{2.20})
 \begin{gather}
X_{1} = \int_{\gamma} \left(\psi_{1}^2 - \psi_{2}^2 \right) dz'
 + \left( \bar{\psi}_{1}^{2} - \bar{\psi}_{2}^2\right) d \bar{z}',
 \nonumber\\
  X_{2} = \int_{\gamma} \left(\psi_{1}^2 + \psi_{2}^2\right)  dz'
 - \left(\bar{\psi}_{1}^2 + \bar{\psi}_{2}^2\right) d \bar{z}',
 \nonumber\\
 X_{3} = - \int_{\gamma} \psi_{1} \psi_{2}  dz'
 + \bar{\psi}_{1} \bar{\psi}_{2} d\bar{z}',
 \label{3.2}
 \end{gather}
 where $\gamma$ is an arbitrary curve, which does not depend
 on the trajectory but only on its endpoints $z$ in ${\mathbb C}$.
 Note that these equations are really written for surfaces
 with $H=1$. To see how to reduce more general cases down to
 this case -- see~\cite{Tai}.
 Note further that   $\psi_1$ and $\psi_2$ are now functions of
 both $z$ and $\bar z$.
 The formulae~(\ref{3.1}) are the starting point of our analysis of
CMC-surfaces in this paper, and according to~\cite{K}, we will
 refer to system~(\ref{3.1}) as the generalised Weierstrass (GW) system.

 In this paper, we examine certain aspects of CMC-surfaces in
 ${\R}^n$ in the context of relating them to solutions of
 low dimensional sigma models. In
 particular, we focus our attention on constructing a Weierstrass
 representation for generic two-dimensional
 surfaces immersed in ${\R}^8$, whose explicit form
 has not been known up to now. For the sake of convenience
 our investigation starts with a derivation of the
 position vector $X (z, \bar{z})$ of a~surface in ${\R}^3$ from
 the Lax pair for a GW system~(\ref{3.1}). As it was shown in~\cite{Bra3}
 the GW system~(\ref{3.1}) is in a one-to-one correspondence
 with the solutions of the equations of the completely integrable
 two-dimensional
 Euclidean $CP^1$ sigma model
 \begin{equation}\label{3.3}
 [ \partial \bar{\partial} P, P ] = 0,
 \end{equation}
 where $P$ is a projector
 \begin{equation}\label{3.4}
 P = \frac{1}{A} \left(
 \begin{array}{cc}
 1  & \bar{w}   \\
 w  & |w|^2     \\
 \end{array}\right)
 ,\qquad
 A = 1 + |w|^2,
 \end{equation}
 or equivalently, the solutions of
 \begin{equation}\label{3.5}
 \partial \bar{\partial} w - \frac{2 \bar{w}}{1 + |w|^{2}} \partial w
 \bar{\partial} w = 0.
 \end{equation}
 In \cite{Ken} it was shown that
 if $\psi_{1}$ and $\psi_{2}$ are solutions of the GW
 system~(\ref{3.1}),
 then function~$w$ defined by
 \begin{equation}\label{3.6}
 w= \frac{\psi_{1}}{\bar{\psi}_{2}},
 \end{equation}
 is a solution of the equations of the $CP^1$ sigma model, namely,~(\ref{3.5}).
The converse is also true~\cite{Bra3}. Thus, if $w$ is a
 solution of~(\ref{3.5}),
 then
 $\psi_{1}$ and $\psi_{2}$ of the GW system (\ref{3.1}) have the form
 \begin{equation}\label{3.7}
 \psi_{1} = \epsilon w \frac{(\bar{\partial} \bar{w})^{1/2}}{1 + |w|^2},
 \qquad
 \psi_{2} = \epsilon \frac{(\partial w)^{1/2}}{1 + |w|^2},
 \qquad
 p = \frac{|\partial w|}{1 + |w|^2},
 \qquad \epsilon = \pm 1.
 \end{equation}
 Note that
 equation (\ref{2.9}) with $H=1$
 for the unit normal vector $N= (n_{1}, n_{2}, n_{3})$ to
 a~CMC-surface adopts the well known form of the equation of
 the $SO(3)$ sigma model
 \begin{equation}\label{3.8}
 \partial \bar{\partial} N + (\partial N, \bar{\partial} N) N = 0,
 \qquad
 (N,N)=1.
 \end{equation}
 Combining the map of the unit vector $N$ onto the unit sphere $S^2$
 with the stereographic projection, we obtain the Gauss map
 \begin{equation}\label{3.9}
 w = \frac{n_{1} + in_{2}}{1 + n_{3}} = \frac{\psi_{1}}{\bar{\psi}_{2}}.
 \end{equation}
 which satisfies the $CP^1$ model equation (\ref{3.5}).

 However, as shown by Zakharov and Mikhailov  \cite{book10} the equation~(\ref{3.3})
 can be considered as a compatibility condition for  two linear
 spectral problems
 \begin{equation}\label{3.10}
 \partial \Psi = \frac{2}{1 + \lambda} [ \partial P, P] \Psi,
 \qquad
 \bar{\partial} \Psi = \frac{2}{1- \lambda} [ \bar{\partial} P, P ],
 \end{equation}
 where $\lambda \in {\mathbb C}$ is a spectral parameter.
 The compatibility condition for (\ref{3.10}) can also be written
 the form of a conservation law
 \begin{equation}\label{3.11}
 \partial K - \bar{\partial} K^{\dagger} = 0,
 \end{equation}
 where the traceless 2 by 2 matrices $K$ and $K^{\dagger}$
 expressed in terms of $w$ have the form
 \begin{gather}
K = [ \bar{\partial} P, P] = \frac{1}{A^2}  \left(
 \begin{array}{cc}
 \bar{w} \bar{\partial} w - w \bar{\partial} \bar{w} &
 \bar{\partial} \bar{w} + \bar{w}^2 \partial w  \vspace{1mm}\\
 - \bar{\partial} w - w^2 \bar{\partial} \bar{w} &
 w \bar{\partial} \bar{w} - \bar{w} \bar{\partial} w  \\
 \end{array}  \right)
 \nonumber\\
 - K^{\dagger} = [ \partial P, P] = \frac{1}{A^2}   \left(
 \begin{array}{cc}
 \bar{w} \partial w - w \partial \bar{w}  &
 \partial \bar{w} + \bar{w}^2 \partial w    \vspace{1mm}\\
 - \partial w - w^2 \partial \bar{w}      &
 w \partial \bar{w} - \bar{w} \partial w      \\
 \end{array}    \right),
 \label{3.12}
 \end{gather}
 and the Hermitian conjugate is denoted by $\dagger$.

 Next we derive the explicit form of matrices $K$ and
 $K^{\dagger}$ in terms of $\psi_{1}$ and $\psi_{2}$
 in order to find the corresponding conservation laws for
 the GW system~(\ref{3.5}). For computational purposes, it is
 useful to express the first derivatives of $w$ in terms
 of $\psi_{1}$ and $\psi_{2}$. Note that in the $CP^1$ case
 the quantity~$J$
 defined in~ (\ref{2.8}) is given by
 $$
 J =\frac{\partial W \partial \bar{W}}{A\sp2}.
 $$
 Using (3.2) we find that it is given by
 \begin{equation}\label{3.13}
 J =\bar{\psi}_{1} \partial
 \psi_{2} - \psi_{2} \partial \bar{\psi}_{1},
 \end{equation}
 and so is a holomorphic function, i.e.\
 which satisfies,
 \begin{equation}\label{3.14}
 \bar{\partial} J = \bar{\partial} (\bar{\psi}_{1}
 \partial \psi_{2} - \psi_{2} \partial \bar{\psi}_{1})
 = - p \partial p + p \partial p = 0,
 \end{equation}
 whenever (\ref{3.1}) holds.
 Actually, in the $CP^1$ case $J$ is a component of the energy-momentum
 tensor.
 Note also that~(3.13) is different from~(2.23).
 It becomes, formally, equal to it under the
 substitution $\bar \psi_1 \rightarrow \psi_1$.
  Using equations (\ref{3.1}), (\ref{3.6}) and
 (\ref{3.13})
 we can express the first derivatives of $w$
 in terms of $\psi_{1}$, $\psi_{2}$ and $J$
 \begin{equation}\label{3.15}
 \partial w = A^{2} \psi_{2}\sp2 ,\qquad
 \bar{\partial} w = - \bar{J} \bar{\psi}_{2}^{-2},
 \end{equation}
 where
 \begin{equation}\label{3.16}
 A= 1 + \frac{|\psi_{1}|^2}{|\psi_{2}|^2}.
 \end{equation}
 As a consequence of (\ref{3.11}) and (\ref{3.12})
 we find that the GW system
 possesses at least three further conservation laws
 \begin{gather}
\partial (\psi_{1} \bar{\psi}_{2} + \bar{R} \bar{\psi}_{1} \psi_{2})
 - \bar{\partial} ( \bar{\psi}_{1} \psi_{2} + R \psi_{1} \bar{\psi}_{2}) =
 0, \nonumber\\
\partial \left(\psi^2_{1} - \bar{R} \psi_{2}^{2}\right)
 + \bar{\partial} \left( \psi_{2}^{2} - R \psi_{1}^{2}\right) = 0,
 \nonumber\\
\partial \left( \bar{\psi}_{2}^{2} - \bar{R} \bar{\psi}_{1}^{2}\right)
 + \bar{\partial} \left( \bar{\psi}_{1}^{2} - R \bar{\psi}_{2}^{2}\right) = 0,
 \label{3.17}
 \end{gather}
 where we have introduced
 \begin{equation}\label{3.18}
 R = \frac{J}{p^2}.
 \end{equation}
 Note that formulae (\ref{3.17}) differ from the conservation laws
 derived in~\cite{K} as they contain additional terms involving~$R$.
 If we put $R=0$ in equations (\ref{3.17}) then we recover the expressions
 given in~\cite{K}
 \begin{equation}\label{3.19}
 \partial (\psi_{1} \bar{\psi}_{2})- \bar{\partial}
 ( \bar{\psi}_{1} \psi_{2}) = 0, \qquad
 \partial \left(\psi_{1}^{2}\right) + \bar{\partial} \left(\psi_{2}^{2}\right) =0,
 \qquad
 \partial \left(\bar{\psi}_{2}^{2}\right) + \bar{\partial} \left( \bar{\psi}_{1}^{2}\right) = 0.
 \end{equation}
 As a result of the conservation laws (\ref{3.17}), we can introduce three
 real-valued functions~$X_{i} (z, \bar{z})$ given by
 \begin{gather}
X_{1} = \frac{i}{2} \int_{\gamma}
 \left[ \bar{\psi}_{1}^{2} + \psi_{2}^{2}  - R \left( \psi_{1}^{2}
 + \bar{\psi}_{2}^{2}\right)\right] dz'
 -\left [ \psi_{1}^{2} + \bar{\psi}_{2}^{2} - R \left(\bar{\psi}_{1}^{2}
 + \psi_{2}^{2}\right) \right]  d \bar{z}', \nonumber\\
X_{2} = \frac{1}{2} \int_{\gamma}
\left[\bar{\psi}_{1}^{2} - \psi_{2}^{2} + R \left(\psi_{1}^{2} -
 \bar{\psi}_{2}^{2}\right)\right]  dz'
 + \left[ \psi_{1}^{2} - \bar{\psi}_{2}^{2} + R \left( \bar{\psi}_{1}^{2} -
 \psi_{2}^{2}\right)\right] d \bar{z}',
 \nonumber\\
X_{3} = - \int_{\gamma} \left[ \bar{\psi}_{1} \psi_{2} + R \psi_{1}
 \bar{\psi}_{2}\right] dz'
 + \left[ \psi_{1} \bar{\psi}_{2} + \bar{R} \bar{\psi}_{1} \psi_{2}\right] d
 \bar{z}',
 \label{3.20}
 \end{gather}
 where $\gamma$ is any curve from a fixed point $z$ in ${\mathbb C}$. The
 functions $X_{i}$, $i=1,2,3$ can be considered as components of a
 position vector of a surface locally parametrised  by $z$ and~$\bar{z}$
 and
 immersed in ${\R}^{3}$
 \begin{equation}\label{3.21}
 X(z, \bar{z}) = ( X_{1} (z, \bar{z}), X_{2} (z, \bar{z}), X_{3} (z,
 \bar{z})).
 \end{equation}
 Using conformal changes of coordinates on the surface ${\mathbb M}^2$ we can,
 without
 loss of generality, put $J=1$ (when, of course $J\ne0$).
  As a consequence it is easy to show that
 representation (\ref{3.20}) with $R=1/p^2$ cannot be reduced
  to the Weierstrass formulae (\ref{3.2}). This means, as we
 will see
 latter, that the additional terms involving $R$ play an important role in
 the construction of surfaces in ${\R}^3$.

 The tangents and the normal unit vector to the surface ${\mathbb M}^2$ are
 given by
 \begin{gather}
 \partial X = \left(i \left[ \bar{\psi}_{1}^{2} + \psi_{2}^{2} - R\left( \psi_{1}^{2} +
 \bar{\psi}_{2}^{2}\right)\right], \left[ \bar{\psi}_{1}^{2} - \psi_{2}^{2} + R
\left(\psi_{1}^{2} -
 \bar{\psi}_{2}^{2}\right)\right],\right.\nonumber\\
\qquad\qquad \left.{} -2 \left(\bar{\psi}_{1} \psi_{2} + R \psi_{1}
 \bar{\psi}_{2}\right)\right),
 \nonumber\\
 \bar{\partial} X = \left(-i\left[\psi_{1}^{2} + \bar{\psi}_{2}^{2} - \bar{R}
\left(\bar{\psi}_{1}^{2} + \psi_{2}^{2}\right)\right], \left[ \psi_{1}^{2} - \bar{\psi}_{2}^{2}
 + \bar{R} \left(\bar{\psi}_{1}^{2} - \psi_{2}^{2}\right)\right],\right.\nonumber\\
\qquad\qquad \left.{} -2 \left(\psi_{1}
 \bar{\psi}_{2} + \bar{R} \bar{\psi}_{1} \psi_{2}\right)\right),
 \label{3.22}
 \end{gather}
 and
 \begin{equation}\label{3.23}
 N = \frac{1}{p} \left( i \left(\bar{\psi}_{1} \bar{\psi}_{2} - \psi_{1} \psi_{2}\right),
 \bar{\psi}_{1} \bar{\psi}_{2} + \psi_{1} \psi_{2},
 |\psi_{1}|^{2} - |\psi_{2}|^{2}\right),
 \end{equation}
 respectively. The first and second fundamental forms of the surface ${\mathbb M}^2$
 are given by
 \begin{gather}
I = (d X, dX) = 4 \left(J dz^2 + p^2\left(1 + |R|^2\right) dz \, d \bar{z} + \bar{J} d
 \bar{z}^{2}\right), \label{3.24}
\\
II = \left(d^2 X, N\right) = \left(4 J + R + \bar{R}\right) dz^2 + \left(2 p + i (R - \bar{R})\right) dz
 d \bar{z}+ \left(4 \bar{J} - R - \bar{R}\right) d \bar{z}^2.
 \nonumber
 \end{gather}
 These quadratic forms contain the Hopf differential $Jdz^2$ and are
 invariant under any conformal changes of coordinates. The Gauss and
 mean curvatures are
 \begin{equation}\label{3.25}
 K = - p^{-2} \partial \bar{\partial} \ln p,
 \qquad H = 1,
 \end{equation}
 respectively.

 Note that  if $J=0$ then $R=0$ and so  the components of
 the
 fundamental forms (3.24) become
 \begin{equation}\label{3.26}
 g_{12} = 2 p^2, \qquad
 g_{11} = g_{22} = 0, \qquad
 b_{12}=2 p, \qquad
 b_{11}= b_{22}=0.
 \end{equation}
 In this case the solutions of GW system (\ref{3.1}) expressed in terms
 of $w$ are represented by~(\ref{3.7}), where $w(z)$ is any holomorphic
 function. According to~\cite{book11}, the energy
 \begin{equation}\label{3.27}
 E = \iint_{\!\!D} \frac{\partial w \bar{\partial} w}{1 + |w|^{2}}\,
 dz \wedge d \bar{z},
 \end{equation}
 is
 finite when the function $w(z)$ is a ratio of polynomials in $z$.
 Geometrically,
 such functions~$\psi_{i}$ parametrise an immersed sphere $S^2 \subset
 {\R}^3$,
 since $J=0$ implies the proportionality of fundamental forms $I$ and
 $II$.

 Note also that if in equations (\ref{3.24}) we put $J \neq 0$ then there
 is
 no
 conformal
 immersion of surfaces in ${\R}^3$. Hence, equations (\ref{3.22}) and
 (\ref{3.23}) imply that the representation (\ref{3.20}) can admit
 different CMC-surfaces which globally have the same Gauss map~(\ref{3.9}).
 This is due to the fact that the tangent vectors ${\partial X}$ and
 ${\bar\partial} X$ depend  on $J$ while the unit normal vector
 $N$ is independent of $J$. In~\cite{Abe,Hof}, using the isometric
 immersions, formulae similar to (\ref{3.22}) and (\ref{3.23}) for
 particular cases of isothermic surfaces have been discussed.

 Let us now discuss the meaning of conservation laws (\ref{3.17}).
 As $J$ is a
 holomorphic function so, according to (\ref{2.11}), we are dealing with
 CMC-surfaces.
 If the $CP^1$ model is defined over $S^2$ then solutions~$w$ of
 (\ref{3.5})
 are either holomorphic or antiholomorphic functions and so $J=0$.
 However, if
 the
 $CP^1$ model is defined on ${\R}^2$ then the function~$w$ is not necessarily
 holomorphic or antiholomorphic and $J \neq 0$.
 Note that when $J\ne0$ the solutions are defined on ${\R}^2/\{a\}$,
 where $\{a\}$ is a small set of points of ${\R}^2$.
 Subtracting~(\ref{3.19}) from (\ref{3.17}) and introducing new
 independent
 variables $\eta$ and $\bar{\eta}$ according to
 \begin{equation}\label{3.28}
 d \eta = J^{1/2} \, dz, \qquad
 d \bar{\eta} = \bar{J}^{1/2} d \bar{z},
 \qquad
 \bar{\partial} J = 0,
 \end{equation}
 we obtain the following set of expressions:
 \begin{gather}
|J|^{2} \left[ \partial_{\eta} \left(\frac{\bar{\psi}_{1} \psi_{2}}{p^2}\right)
 - \bar{\partial}_{\bar{\eta}} \left(\frac{\psi_{1} \bar{\psi}_{2}}{p^2}\right)\right]  =
 0,
 \qquad
|J|^{2} \left[ \partial_{\eta} \left(\frac{\psi_{2}^{2}}{p^2}\right)
 + \bar{\partial}_{\bar{\eta}} \left( \frac{\psi_{1}}{p^2}\right)\right] = 0,
 \nonumber\\
|J|^{2} \left[ \partial_{\eta} \left( \frac{\bar{\psi}_{1}^{2}}{p^2}\right)
 + \bar{\partial}_{\bar{\eta}} \left( \frac{\bar{\psi}_{2}^{2}}{p^2}\right)\right] = 0,
 \label{3.29}
 \end{gather}
 where the derivatives are abbreviated $\partial_{\eta}= \partial/
 \partial \eta$
 and $\bar{\partial} \bar{\eta} = \partial / \partial
 \bar{\eta}$.
  Equations (\ref{3.29}) suggest that we should consider two
 separate
 cases, namely
 $J=0$ which has been already treated in~\cite{K}
 and $J \neq 0$. In the latter case,
 under the change of variables (\ref{3.28}) the GW system~(\ref{3.1})
 adopts
 the form
 \begin{equation}\label{3.30}
 \partial_{\eta} \psi_{1} = \frac{p}{J}\, \psi_{2},
 \qquad
 \bar{\partial}_{\bar{\eta}} \psi_{2} = - \frac{p}{\bar{J}}\, \psi_{1},
 \end{equation}
 and the expression for $J$, given in (\ref{3.13}), provides
 the following differential constraint~(DC) on
 $\psi_{1}$ and $\psi_{2}$
 \begin{equation}\label{3.31}
 \bar{\psi}_{1} \partial_{\eta} \psi_{2} - \psi_{2}
 \partial_{\eta} \bar{\psi}_{1} = 1.
 \end{equation}

 Keeping in mind that the complex coordinates
 $\eta$ and $\bar{\eta}$ are is defined up to a conformal
 transformation, we can without loss of generality put $J=1$.
 If $\psi_{1} \neq 0$ then the system~(\ref{3.30}), subject to DC
 (\ref{3.31}), can be written in an equivalent form
 \begin{equation}\label{3.32}
 \partial \psi_{1} = p \psi_{2}, \qquad
 \partial \psi_{2} = \bar{\psi}^{-1} ( 1 + \psi_{2} \partial
 \bar{\psi}_{1}), \qquad
 \bar{\partial} \psi_{2} = - p \psi_{1}.
 \end{equation}
 The compatibility condition for (\ref{3.32}) does not imply any new DC on
 first order
 derivatives of $\psi_{1}$. Hence, the system (\ref{3.32}) is integrable
 and
 the derivatives
 $\bar{\partial} \psi_{1}$ and $\partial \psi_{1}$ are undetermined.

 The Gaussian curvature and mean curvature are
 \begin{equation}\label{3.33}
 K = \left(|\psi_{2}|^2 - |\psi_{1}|^{2}\right)
\left [ |\bar{\partial} \psi_{1}|^{2} + |\psi_{2}|^{-2}
 \left(1 + \psi_{2} \partial \bar{\psi}_{1} +\bar{\psi}_{2} \bar{\partial}
 \psi_{1}\right)\right], \qquad H=1,
 \end{equation}
 respectively.

 Note that the equations  of the complex frame (\ref{2.6}) are specified
 by DC~(\ref{3.31}) and are
 compatible with the scalar products (\ref{2.3}) and (\ref{2.4}).
 After the change of dependent variables
 \begin{equation}\label{3.34}
 p = e^{\varphi/2},
 \end{equation}
 the corresponding Gauss--Codazzi equations (\ref{2.10})--(\ref{2.11})
 take the
 form of the elliptic Sh-Gordon equation
 \begin{equation}\label{3.35}
 \partial \bar{\partial} \varphi + 4 \sinh \varphi = 0.
 \end{equation}
 Hence the CMC-surfaces are determined by formulae (\ref{3.20}), where
  $\psi_{1}$ and $\psi_{2}$ have to obey
 equations (\ref{3.32}) with $p$ determined by (\ref{3.34})
 and (\ref{3.35}).
 In terms of arbitrary conformal coordinates, we have proved that
  $(\psi_{1}, \psi_{2}, p )$ can be viewed as the Weierstrass
 data
 of the CMC-surface ${\mathbb M}^2$ in ${\R}^3$.

 To summarize: the generalised
 Weierstrass representation for the
 immersion of a CMC-surface into ${\R}^3$ is described
 by formulae (\ref{3.20}), where~$\psi_{1}$ and $\psi_{2}$
 obey the GW system of equations~(\ref{3.1}).

 Let us add also, as shown in~\cite{Fer}, that under the changes of
 independent variables~(\ref{3.28})
 and dependent variables $S= 2 p^2 |J|^{-1}$, the GW system (\ref{3.1}) is
 decoupled
 into a direct sum of  equations
 \begin{equation}\label{3.36}
 \partial_{\eta} \bar{\partial}_{\bar{\eta}} \ln S = S^{-1} - S,
 \qquad
 \bar{\partial}_{\bar{\eta}} J = 0,
 \qquad
 \partial_{\eta} \bar{J} = 0.
 \end{equation}
 Hence, the GW system (\ref{3.1}) is completely integrable. Examples of
 $N$
 solitons solutions can be found in~\cite{Bra3}.

 \section{The $\boldsymbol{CP^2}$ maps and the  Weierstrass representation\\
 for surfaces in  eight  dimensional Euclidean spaces}

 The aim of this section is to demonstrate a connection
 between the recently proposed generalised Weierstrass (GW) system~\cite{G}
 \begin{gather}
\partial \psi_{1} = \left(1+ \frac{|\psi_{2}|^2}{|\varphi_{2}|^2}\right)
 \bar{\varphi}_{1} \bar{Q} - \frac{1}{2} \left[
 \frac{\psi_{1} \bar{\psi}_{2}}{\varphi_{2}}
 + \frac{|\psi_{1}|^2}{|\varphi_{1}|^2}
 \frac{\bar{\varphi}_{2}^2}{\bar{\varphi}_{1}}\right] \bar{P},
 \label{4.1}\\
\partial \psi_{2} = \left(1 + \frac{|\psi _{1}|^2}{|\varphi_{1}|^2}\right)
 \bar{\varphi}_{2} \bar{P} - \frac{1}{2} \left[
 \frac{\bar{\psi}_{1} \psi_{2}}{\varphi_{1}}
 + \frac{|\psi_{2}|^2}{|\varphi_{2}|^2}
 \frac{\bar{\varphi}_{1}^2}{\bar{\varphi}_{2}}\right] \bar{Q},
 \label{4.2}\\
\bar{\partial} \varphi_{1} =
 - \frac{1}{2} \left[\left( \frac{\psi_{2}}{\bar{\varphi}_{2}} P
 + 2 \frac{\psi_{1}}{\bar{\varphi}_{1}} Q\right) \varphi_{1} + P \psi_{1}
 \frac{\varphi_{2}^2}{|\varphi_{1}|^2}\right],
 \label{4.3}\\
\bar{\partial} \varphi_{2} = - \frac{1}{2} \left[\left(
 \frac{\psi_{1}}{\bar{\varphi}_{1}} Q
 + 2 \frac{\psi_{2}}{\bar{\varphi}_{2}} P\right) \varphi_{2} + Q \psi_{2}
 \frac{\varphi_{1}^2}{|\varphi_{2}|^2}\right],
 \label{4.4}
 \end{gather}
 where
 \begin{equation}
P = \frac{\psi_{1} \bar{\psi}_{2} \bar{\varphi}_{1}}{\varphi_{2}}
 + \left( |\varphi_{2}|^2 + |\psi_{2}|^2\right) \frac{\bar{\varphi}_{2}}{\varphi_{2}},
 \qquad
Q = \frac{\bar{\psi}_{1} \psi_{2} \bar{\varphi}_{2}}{\varphi_{1}}
 + \left(|\varphi_{1}|^2 + |\psi_{1}|^2 \right) \frac{\bar{\varphi}_{1}}{\varphi_{1}}
 \label{4.5}
 \end{equation}
 and the equations of the $CP^2$ sigma model \cite{book11}
 \begin{gather}
\partial \bar{\partial} w_{1} - \frac{2 \bar{w_{1}}}{A}
 \partial w_{1}
 \bar{\partial} w_{1} - \frac{\bar{w}_{2}}{A}
 (\partial w_{1} \bar{\partial} w_{2} + \bar{\partial} w_{1} \partial
 w_{2}) = 0,
 \label{4.6}\\
\partial \bar{\partial} w_{2} - \frac{2 \bar{w}_{2}}{A} \partial w_{2}
 \bar{\partial} w_{2} - \frac{\bar{w}_{1}}{A} (\partial w_{1}
 \bar{\partial} w_{2}
 + \bar{\partial} w_{1} \partial w_{2}) = 0,
 \label{4.7}\\
A = 1 + |w_{1}|^2 + |w_{2}|^2.
 \label{4.8}
 \end{gather}

 Next, we exploit this connection and use the conservation laws
 for the GW system
 (\ref{4.1})--(\ref{4.4})
 to define real valued functions $X^{i} (z, \bar{z})$, $i=1, \ldots, 8$
 in terms of functions $\varphi_{\alpha}$, $\psi_{\alpha}$, $\alpha=1,2$
 which
 are identified as the coordinates in $8$-dim Euclidean space ${\R}^8$.
 The formulae (\ref{4.1})--(\ref{4.4})
 and (\ref{4.6})--(\ref{4.8}) are the starting point for our
 analysis.
 In this paper, when we refer to system (\ref{4.1})--(\ref{4.4}), we
 will describe it as the
 modified  version of the original Weierstrass system~(\ref{3.1}).

 Note that in equations (\ref{4.1})--(\ref{4.4}) only four out of eight
 derivatives of functions $\psi_{i}$ and~$\varphi_{i}$ are known in terms
 of complex functions $\psi_{i}$ and $\varphi_{i}$ and their complex
 conjugates while the others are unspecified. Note also that if the
 functions $\psi_{\alpha}$ tend to $\psi/\sqrt{2}$ and $\varphi_{\alpha}$
 tend to $\varphi$, i.e.\ then the system (\ref{4.1})--(\ref{4.4})
  reduces to the
 Weierstrass formulae (\ref{3.1}) for the CMC-surfaces immersed
 in ${\R}^3$
 \begin{gather*}
\partial \psi = \left(|\psi|^2 + |\varphi|^2\right) \varphi,
 \qquad
\bar{\partial} \varphi = - \left(|\psi|^2 + |\varphi|^2\right) \psi.
 \end{gather*}
 In terms of $w_{i}$, $i=1,2$
 the above limit takes the form
 \begin{equation}\label{4.9}
 w_{i} \rightarrow \frac{1}{\sqrt{2}}\, w,
 \qquad i=1,2,
 \end{equation}
 and then the $CP^2$ sigma model (\ref{4.6})--(\ref{4.8}) reduces to the
 $CP^1$ sigma model (\ref{3.5}). These limits characterise the properties
 of
 the solutions of
 systems (\ref{4.1})--(\ref{4.4}) and (\ref{4.6})--(\ref{4.8}).

 First we show that there exists a one to one correspondence between
 $GW$ system (\ref{4.1})--(\ref{4.4}) and the equations of the
 $CP^2$ sigma model (\ref{4.6})--(\ref{4.8}). For this purpose, we define
 two  new complex valued functions
 \begin{equation}\label{4.10}
 w_{1} = \frac{\psi_{1}}{\bar{\varphi}_{1}},
 \qquad
 w_{2} = \frac{\psi_{2}}{\bar{\varphi}_{2}}
 \end{equation}
 and using the GW system (\ref{4.1})--(\ref{4.4}), we obtain
 \begin{gather}
\partial w_{1} = A \left[ w_{1} \bar{w}_{2} \varphi_{2}^2
 + \left( 1 + |w_{1}|^2\right) \varphi_{1}^2\right],
 \nonumber\\
\partial w_{2} = A \left[ \bar{w}_{1} w_{2} \varphi_{1}^2
 + \left(1 + |w_{2}|^2\right) \varphi_{2}^2\right].
 \label{4.11}
 \end{gather}
 These relations generate the following transformation from the
 variables $(w_{1}, w_{2})$ and their derivatives to the
 variables $( \varphi_{1}, \varphi_{2}, \psi_{1}, \psi_{2})$
 \begin{gather}
\varphi_{1} = \epsilon A^{-1}
 \left[ \left( 1 + |w_{2}|^2 \right) \partial w_{1} - w_{1} \bar{w}_{2}
 \partial w_{2}\right]^{1/2},
 \label{4.12}\\
\varphi_{2} = \epsilon A^{-1} \left[- \bar{w}_{1} w_{2} \partial w_{1}
 + \left( 1 + |w_{1}|^2\right) \partial w_{2} \right] ^{1/2},
 \label{4.13}\\
\psi_{1} = \epsilon w_{1} A^{-1}
 \left[\left(1 + |w_{2}|^2\right) \bar{\partial} \bar{w}_{1} - \bar{w}_{1} w_{2}
 \bar{\partial} \bar{w}_{2} \right]^{1/2},
 \label{4.14}\\
\psi_{2} = \epsilon w_{2} A^{-1} \left[- w_{1} \bar{w}_{2}
 \bar{\partial} \bar{w}_{1} + \left(1 + |w_{1}|^2\right)
 \bar{\partial} \bar{w}_{2}\right]^{1/2}.
 \label{4.15}
 \end{gather}

 We can now state the following: if the complex valued functions
 $(\varphi_{1}, \varphi_{2}, \psi_{1}, \psi_{2})$ are
 solutions of GW system (\ref{4.1})--(\ref{4.4}), then the  functions
 $(w_{1}, w_{2})$, defined by (\ref{4.10}), solve
 the equations of the $CP^2$ sigma model (\ref{4.6})--(\ref{4.8}).

 Conversely, if the complex valued functions $(w_{1}, w_{2})$ are
 solutions of the $CP^2$ sigma model equations (\ref{4.6})--(\ref{4.8}),
 then
 the
 complex valued functions $(\varphi_{1}, \varphi_{2}, \psi_{1}, \psi_{2})$
 defined by (\ref{4.12})--(\ref{4.15}) in
 terms of functions $(w_{1}, w_{2})$ and their 1st
 derivatives satisfy the GW system (\ref{4.1})--(\ref{4.4}).

 The proof of our statement is straightforward.  The differentiation of
 equations (\ref{4.11}) with respect to $z$
 and $\bar{z}$, respectively, yields
 \begin{gather*}
\partial \bar{\partial} w_{1} = A \left[ \bar{w}_{2} \varphi_{2}^2
 \bar{\partial}
 w_{1} + w_{1} \varphi_{2}^2 \bar{\partial}
 \bar{w}_{2} + 2 w_{1} \bar{w}_{2} \varphi_{2}
 \bar{\partial} \varphi_{2}
 +2 \left( 1 + |w_{1}|^2\right) \varphi_{1} \bar{\partial} \varphi_{1}\right.
\\
\left. \phantom{\partial \bar{\partial} w_{1} =}{} +
\left( \bar{w}_{1} \bar{\partial} w_{1} + w_{1} \bar{\partial}
 \bar{w}_{1}\right)
 \varphi_{1}^2 \right]
 + \left[ w_{1} \bar{w}_{2} \varphi_{2}^2 + \left(1 + |w_{1}|^2\right) \varphi_{1}^2\right]\\
\phantom{\partial \bar{\partial} w_{1} =}{}\times
 \left( \bar{w}_{1} \bar{\partial} w_{1} + w_{1} \bar{\partial} \bar{w}_{1}
 + \bar{w}_{2} \bar{\partial} w_{2} + w_{2} \bar{\partial} \bar{w}_{2}\right),
 \end{gather*}
 and
 \begin{gather*}
\partial \bar{\partial} w_{2} = A \left[ w_{2} \varphi_{1}^2 \bar{\partial}
 \bar{w}_{1}
 + \bar{w}_{1} \varphi_{1}^2 \bar{\partial} w_{2} + 2 \bar{w}_{1} w_{2}
 \varphi_{1}
 \bar{\partial} \varphi_{1}
 + 2 \left(1 + |w_{2}|^2\right) \varphi_{2} \bar{\partial} \varphi_{2}\right.
\\
\left.\phantom{\partial \bar{\partial} w_{2} =}{}+
\left( \bar{w}_{2} \bar{\partial} w_{2} + w_{2} \bar{\partial}
 \bar{w}_{2}\right)
 \varphi_{2}^2\right]
 + \left[ \bar{w}_{1} w_{2} \varphi_{1}^{2} + \left(1 + |w_{2}|^2 \right) \varphi_{2}^2\right]\\
 \phantom{\partial \bar{\partial} w_{2} =}{}\times
\left( \bar{w}_{1} \bar{\partial} w_{1} + w_{1} \bar{\partial} \bar{w}_{1} +
 \bar{w}_{2}
 \bar{\partial} w_{2} + w_{2} \bar{\partial} \bar{w}_{2}\right)
 \end{gather*}
 and their respective complex conjugate equations.

 Substituting (\ref{4.11}) and
 into the left-hand side of the first
 equation (\ref{4.6}), we obtain
 \begin{gather}
\partial \bar{\partial} w_{1} - \frac{2 \bar{w}_{1}}{A} \partial w_{1}
 \bar{\partial} w_{1} - \frac{\bar{w}_{2}}{A} (\partial w_{1}
 \bar{\partial} w_{2}
 + \bar{\partial} w_{1} \partial w_{2})
 \nonumber\\
\qquad {}= \left[ A w_{1} \varphi_{2}^2 + w_{1} |w_{2}|^2 \varphi_{2}^2 + \left(1 +
 |w_{1}|^2\right) w_{2}
 \varphi_{1}^2\right]
\left( w_{1} \bar{w}_{2} \bar{\varphi}_{1}^2 +\left(1 + |w_{2}|^2 \right)
 \bar{\varphi}_{2}^2\right)
 \nonumber\\
\qquad{}+ \left[ A w_{1} \varphi_{1}^2 + w_{1}^2 \bar{w}_{2}^2 \varphi_{2}^2 +
\left( 1 + |w_{1}|^2\right) w_{1} \varphi_{1}^2\right] \left( \bar{w}_{1} w_{2}
 \bar{\varphi}_{2}^2
 + \left( 1 + |w_{1}|^2 \right) \bar{\varphi}_{1}^2\right)
 \nonumber\\
\qquad {} + 2 \left[ w_{1} \bar{w}_{2} \varphi_{2} \bar{\partial} \varphi_{2}
 + \left( 1 + |w_{1}|^2 \right) \varphi_{1} \bar{\partial} \varphi_{1}\right].
 \label{4.16}
 \end{gather}

 Making use of the equations (\ref{4.1})--(\ref{4.4}) we find that
 the equation (\ref{4.16})
 is satisfied identically.  An analogous
 result takes place for the second equation (\ref{4.7}), since
 the $CP^2$ sigma
 model equations (\ref{4.6})--(\ref{4.8}) are invariant under
\begin{equation}
 \label{4.17}
 w_{1} \leftrightarrow w_{2}
 \end{equation}
 This observation implies that the left-hand side of (\ref{4.7})
 vanishes as well
 whenever (\ref{4.1})--(\ref{4.4}) holds.

 Conversely, differentiating (\ref{4.12}) with respect to $\bar{z}$ and
 using~(\ref{4.12}),
 we get
 \begin{gather}
\bar{\partial} \varphi_{1} = \frac{1}{2 A^2 \varphi_{1}}
 \left[ \left( \bar{w}_{2} \bar{\partial} w_{2} + w_{2} \bar{\partial} \bar{w}_{2}\right)
 \partial w_{1} + \left(1 + |w_{2}|^2\right) \partial \bar{\partial} w_{1}\right.\nonumber\\
\left.\phantom{\bar{\partial} \varphi_{1} =}{}
-\left(\bar{w}_{2} \bar{\partial} w_{1} + w_{1} \bar{\partial}
 \bar{w}_{2}\right) \partial w_{2}
 - w_{1} \bar{w}_{2} \partial \bar{\partial} w_{2}\right]
 \nonumber\\
\phantom{\bar{\partial} \varphi_{1} =}{}- \frac{\varphi_{1}}{A}
\left[\bar{w}_{1} \bar{\partial} w_{1} + w_{1} \bar{\partial}  \bar{w}_{1}
 + \bar{w}_{2} \bar{\partial} w_{2} + w_{2} \bar{\partial} \bar{w}_{2}\right].
 \label{4.18}
 \end{gather}
 Using equations (\ref{4.11}) and (\ref{4.6})--(\ref{4.8}), we can
 eliminate first and second
 derivatives
 of $w_{1}$ and $w_{2}$ in expression (\ref{4.18}) and obtain
 \begin{gather}
\bar{\partial} \varphi_{1} = \frac{1}{2 \varphi_{1}} \left\{ w_{1}^2
 \bar{w}_{2}
 |w_{2}|^2 \bar{\varphi}_{1}^2 \varphi_{2}^2 + \left(1 + |w_{2}|^2\right) w_{1}
 |w_{2}|^2
 |\varphi_{2}|^4 + \left(1 + |w_{1}|^2\right) w_{1} |w_{2}|^2 |\varphi_{1}|^4\right.
 \nonumber\\
\phantom{\bar{\partial} \varphi_{1} =}{}
+\left(1 + |w_{1}|^2 \right)\left( 1 + |w_{2}|^2\right) w_{2} \varphi_{1}^2
 \bar{\varphi}_{2}^2
 - w_{1} | w_{1}|^2 | w_{2}|^2 |\varphi_{1}|^4
 \nonumber\\
\phantom{\bar{\partial} \varphi_{1} =}{}
- \left(1 + |w_{2}|^2\right) |w_{1}|^2
 w_{2} \varphi_{1}^2 \bar{\varphi}_{2}^2
- \left(1 + |w_{2}|^2\right) w_{1}^2 \bar{w}_{2} \bar{\varphi}_{1}^2 \varphi_{2}
 - \left( 1 + |w_{2}|^2\right)^2 w_{1} |\varphi_{2}|^4
 \nonumber\\
\phantom{\bar{\partial} \varphi_{1} =}{}
- 2 A \varphi_{1}^2
\left(w_{2} \bar{\varphi}_{2}^2 + w_{1} \bar{\varphi}_{1}^2\right)
+ \frac{1}{A^2} \left(A \left[ w_{1} \bar{w}_{2}^2 \varphi_{2}^2 + \bar{w}_{2}
\left(1 + |w_{1}|^2 \right) \varphi_{1}^2\right] \right.
 \nonumber\\
\phantom{\bar{\partial} \varphi_{1} =}{}
+ \left(1 + |w_{2}|^2\right)
\left[ w_{1} \bar{w}_{2}^2 \varphi_{2}^2 + \bar{w}_{2} \left(1 + |w_{1}|^2\right)
 \varphi_{1}^2\right]
 \nonumber\\
\left.\phantom{\bar{\partial} \varphi_{1} =}{}
- \left(A + |w_{2}|^2\right) |w_{1}|^2 \bar{w}_{2} \varphi_{1}^2
 - \left(A +1 + |w_{2}|^2\right) w_{1} \bar{w}_{2}^2 \varphi_{2}^2 -2 A \bar{w}_{2}
 \varphi_{1}^2\right) \bar{\partial} w_{2}
 \nonumber\\
\phantom{\bar{\partial} \varphi_{1} =}{}
+ \frac{1}{A^2} \left(\left(1 + |w_{2}|^2 \right) \left[ \left(A + |w_{1}|^2\right) \bar{w}_{2}
 \varphi_{2}^2
 + \left(A + 1 + |w_{1}|^2\right) \bar{w}_{1} \varphi_{1}^2\right]\right.
 \nonumber\\
\phantom{\bar{\partial} \varphi_{1} =}{}
- A \left[ \bar{w}_{1} | w_{2}|^2 \varphi_{1}^2 + \left(1 + |w_{2}|^2\right)
 \bar{w}_{2}
 \varphi_{2}^2\right]
 \nonumber\\
\left.\left.\phantom{\bar{\partial} \varphi_{1} =}{}
-\left[ \bar{w}_{1} | w_{1}|^2 |w_{2}|^2
 \varphi_{1}^2
+ \left( 1 + |w_{2}|^2\right) |w_{1}|^2 \bar{w}_{2} \varphi_{2}^2\right]
 -2 A \bar{w}_{1} \varphi_{1}^2 \right) \bar{\partial} w_{1} \right\}.
 \label{4.19}
 \end{gather}
 Collecting all the coefficients of the derivatives $\bar{\partial} w_{1}$
 and $\bar{\partial} w_{2}$ in expression (\ref{4.19}) we find that
 these coefficients vanish identically. In fact, we have
 \begin{gather*}
\bar{\partial} w_{1} : \ \left(A + 1 + |w_{1}|^2 + \left(A + 1 + |w_{1}|^2\right)
 |w_{2}|^2 - A |w_{2}|^2 - |w_{1}|^2 |w_{2}|^2 - 2 A\right) \bar{w}_{1}
 \varphi_{1}^2\\
\qquad {}+ \left(A + |w_{1}|^2 + \left(A + |w_{1}|^2\right) |w_{2}|^2 - A \left(1 + |w_{2}|^2 \right)
 - \left(1 + |w_{2}|^2\right) |w_{1}|^2 \right) \bar{w}_{2} \varphi_{2}^{2} \equiv 0.\!
 \end{gather*}
 and
 \begin{gather}
\bar{\partial} w_{2} : \ \left( A \left(1 + |w_{1}|^2 \right) + A + |w_{1}|^2 |w_{2}|^2
 - \left(A + |w_{2}|^2\right) |w_{1}|^2  - 2 A\right) \bar{w}_{2} \varphi_{1}^2
 \nonumber\\
\qquad {}+\left( A + 1 + |w_{2}|^2 - A -1 - |w_{2}|^2\right) w_{1} \bar{w}_{2} \varphi_{2}^2
 \equiv 0.
 \label{4.20}
 \end{gather}
 Hence,  (\ref{4.20}) becomes
 \begin{gather*}
\bar{\partial} \varphi_{1} = - \frac{1}{2} \left\{ \bar{w}_{2} w_{1}^2
 \frac{\bar{\varphi}_{1}^2 \varphi_{2}^2}{\varphi_{1}} +
\left(1 + |w_{2}|^2\right) w_{1} \frac{|\varphi_{2}|^4}{\varphi_{1}}
 + 2 A w_{2} \varphi_{1} \bar{\varphi}_{2}^2\right.\nonumber\\
\left.\phantom{\bar{\partial} \varphi_{1} =}{}
 + 2 A w_{1} | \varphi_{1}|^2 \varphi_{1}
- w_{1} | w_{2}|^2 \frac{|\varphi_{1}|^4}{\varphi_{1}}
 -\left( 1 + |w_{2}|^2\right) w_{2} \varphi_{1} \bar{\varphi}_{2} \right\}
 \end{gather*}
 Performing the transformation (\ref{4.10})
 we obtain the first equation of (\ref{4.1}), i.e.\
 \begin{gather}
 \bar{\partial} \varphi_{1} = - \frac{1}{2} \Big\{
 \frac{\varphi_{2}}{\varphi_{1}}
 \bar{\psi}_{2} \psi_{1}^2 + \left( 1 + |w_{2}|^2\right)
 \frac{|\varphi_{2}|^4}{|\varphi_{1}|^2}\nonumber\\
\phantom{\bar{\partial} \varphi_{1} =}{} +
\left( A +1 + |w_{1}|^2\right) | \varphi_{1}|^2 \psi_{1} + \left(A + |w_{1}|^2\right)
 \varphi_{1}
 \bar{\varphi}_{2} \psi_{2} \Big\}.\label{4.21}
 \end{gather}
 Since the equations  (\ref{4.1})--(\ref{4.4}) are invariant under
 \begin{equation}\label{4.22}
 \varphi_{1} \leftrightarrow \varphi_{2},
 \end{equation}
 an analogous result holds for~(\ref{4.2}).
 Differentiation of (\ref{4.10}) with respect to~$z$ gives
 \begin{equation}\label{4.23}
 \partial \psi_{1} = \bar{\varphi}_{1} \partial w_{1} + w_{1}
 \partial \bar{\varphi}_{1}.
 \end{equation}
 Substituting (\ref{4.11}) and the complex conjugate equation of~(\ref{4.1})
 into (\ref{4.23}) we get~(\ref{4.3}).
 Making use of the discrete symmetry (\ref{4.22})
 in (\ref{4.3}), we obtain equation (\ref{4.4}),
 which completes the proof.

 An interesting property of the GW system (\ref{4.1})--(\ref{4.4}) in the
 context of
 the $CP^2$ sigma model (\ref{4.6})--(\ref{4.8}) is the existence of a gauge
 freedom
 in the definition of the variables~$w_{1}$ and $w_{2}$ given by
formula~(\ref{4.10}). This is due to the fact that the numerator and the
 denominator of~(\ref{4.10})
 can be multiplied by any complex functions $f_{i} : {\mathbb C} \rightarrow {\mathbb C}$,
 $i=1,2$.
 This means that if we introduce a new set of complex valued functions
 $(\alpha_{1}, \alpha_{2}, \beta_{1}, \beta_{2})$ which are related to
 functions $(\varphi_{1}, \varphi_{2}, \psi_{1}, \psi_{2})$ in the
 following way
 \begin{equation}\label{4.24}
 \varphi_{i} = f_{i} (z, \bar{z}) \alpha_{i},  \qquad
 \psi_{i} = \bar{f} ( z, \bar{z}) \beta_{i},  \qquad  i=1,2,
 \end{equation}
 then the transformation (\ref{4.24}) leaves the functions $w_{1}$, $w_{2}$
 invariant
 \begin{equation}\label{4.25}
 w_{1} = \frac{\beta_{1}}{\bar{\alpha}_{1}},
 \qquad
 w_{2} = \frac{\beta_{2}}{\bar{\alpha}_{2}}.
 \end{equation}

 We show now that if the complex valued functions $w_{1}$, $w_{2}$ are
 solutions of the $CP^2$ sigma model equations (\ref{4.6})--(\ref{4.8}), then for
 any two
 holomorphic functions $f_{i}$, $i=1,2$ the complex functions
 $(\alpha_{1}, \alpha_{2}, \beta_{1}, \beta_{2})$ defined by
\begin{gather}
\alpha_{1} = \epsilon f_{1}^{-1} A^{-1}
\left[\left( 1 + |w_{2}|^2\right) \partial w_{1} - w_{1} \bar{w}_{2} \partial
 w_{2}\right]^{1/2},
 \nonumber\\
\alpha_{2} = \epsilon f_{2}^{-1} A^{-1}
\left[ - \bar{w}_{1} w_{2} \partial w_{1} + \left( 1 + |w_{1}|^2\right)
 \partial w_{2}\right]^{1/2},
 \nonumber\\
\beta_{1} = \epsilon w_{1} \bar{f}_{1}^{-1} A^{-1}
 \left[ \left( 1 + |w_{2}|^2\right) \bar{\partial} \bar{w}_{1} - \bar{w}_{1} w_{2}
 \bar{\partial} \bar{w}_{2}\right]^{1/2},
 \nonumber\\
\beta_{2} = \epsilon w_{2} \bar{f}_{2}^{-1} A^{-1}
\left [-w_{1} \bar{w}_{2} \bar{\partial} \bar{w}_{1}
 + \left(1 + |w_{1}|^2\right) \bar{\partial} \bar{w}_{2}\right]^{1/2},
 \qquad \bar{\partial}  f_{i}  = 0,
 \label{4.26}
 \end{gather}
 satisfy the GW system (\ref{4.1})--(\ref{4.4}).

 Indeed, the result is obtained directly by substituting
 (\ref{4.24}) and (\ref{4.25})
 into $CP^2$ sigma model equations
 (\ref{4.6})--(\ref{4.8}).
 This leads to differential constraints for the functions $f_{i}$
 and their first derivatives
 \begin{equation}\label{4.27}
\left(f_{i}^2 - 1\right) \bar{\partial} f_{j} = 0,  \qquad
\left(\bar{f}_{i}^2 - 1\right) \partial \bar{f}_{j} = 0,  \qquad  i,j = 1,2.
 \end{equation}
 Hence, the general solutions of this system are given by any
 holomorphic functions $f_{i}$ i.e.\
 \begin{equation}\label{4.28}
 \bar{\partial} f_{i} = 0, \qquad i=1,2
 \end{equation}
  Then invoking the main result of the last section
 we note  that the transformation
 (\ref{4.12})--(\ref{4.15}) becomes the one given by (\ref{4.26}).

 Another interesting property in the context of the $CP^2$ sigma
 model (\ref{4.6})--(\ref{4.8}) and the GW system (\ref{4.1})--(\ref{4.4}) is
 the
 existence of the quantity $J$~\cite{book11} (which is a generalization
 of~(3.17)).
 \begin{equation}\label{4.29}
 J = A^{-2} \left\{ \partial w_{1} \partial \bar{w}_{1} + \partial w_{2}
 \partial \bar{w}_{2} + (\bar{w}_{1} \partial \bar{w}_{2}
 - \bar{w}_{2} \partial \bar{w}_{1}) ( w_{1} \partial w_{2}
 - w_{2} \partial w_{1}) \right\},
 \end{equation}
 whose derivative with respect to $z$
 vanishes identically whenever equations
 (\ref{4.6})--(\ref{4.8}) are
 satisfied
 \begin{equation}\label{4.30}
 \bar{\partial} J = 0.
 \end{equation}
 This means that  $J$, given by (\ref{4.29}),
 is  holomorphic.

 Note that if the functions $(\varphi_{1}, \varphi_{2}, \psi_{1},
 \psi_{2})$ are
 solutions of GW system (\ref{4.1})--(\ref{4.4}), then~$J$, when written in
 terms
 of
 functions $(\varphi_{1}, \varphi_{2}, \psi_{1}, \psi_{2})$, takes the
 form
 \begin{equation}\label{4.31}
 J = \varphi_{1} \partial \bar{\psi}_{1} - \bar{\psi}_{1} \partial
 \varphi_{1}
 + \varphi_{2} \partial \bar{\psi}_{2} - \bar{\psi}_{2} \partial
 \varphi_{2}
 \end{equation}
 and it satisfies
 \begin{equation}\label{4.32}
 \bar{\partial} J = 0,
 \end{equation}
 whenever equations (\ref{4.1})--(\ref{4.4}) hold.

 Next we exploit the observation
 \cite{book10} that the equations of the $CP^2$ sigma model
  (\ref{4.6})--(\ref{4.8}) can be written as the compatibility condition for
 two
 linear spectral problems
 \begin{equation}\label{4.36}
 \partial \Phi = \frac{2}{1 + \lambda} [ \partial P, P] \Phi,
 \qquad
 \bar{\partial} \Phi = \frac{2}{1 - \lambda} [ \bar{\partial} P, P] \Phi,
 \qquad
 \lambda \in {\mathbb C},
 \end{equation}
 where the $3 \times 3$ $P$ is  given by
 \begin{equation}\label{4.37}
 P = A^{-1} M, \qquad   M = \left(
 \begin{array}{ccc}
 1 & w_{1}  & w_{2}   \\
 \bar{w}_{1} &  |w_{1}|^2 &  \bar{w}_{1} w_{2}   \\
 \bar{w}_{2}  &  w_{1} \bar{w}_{2} & |w_{2}|^2   \\
 \end{array}   \right),
 \end{equation}
 and $\lambda$ represents the spectral parameter.
 Using matrix $P$, the compatibility conditions of equations (\ref{4.36}) imply
 \begin{equation}\label{4.38}
 [\partial \bar{\partial} P, P] = 0,
 \end{equation}
 which are satisfied
 whenever equations (\ref{4.6})--(\ref{4.8}) hold. Equivalently,
 formula (\ref{4.38})
 can be rewritten, in a divergent form, as
 \begin{equation}\label{4.39}
 \partial [ \bar{\partial} P, P] + \bar{\partial}
 [\partial P, P ] = 0.
 \end{equation}
 Hence, from equations (\ref{4.37}) and (\ref{4.39})
 we obtain the explicit form
 of the local conservation laws for the $CP^2$ sigma model
 \begin{equation}\label{4.40}
 \partial K + \bar{\partial} L = 0,
 \end{equation}
 where we have introduced the following notation for the traceless
 matrices
 $K$ and $L$:
 \begin{equation}\label{4.41}
 K = \frac{1}{A^2}\, [ \bar{\partial} M, M],  \qquad
 L = -K^\dagger = \frac{1}{A^2}\, [ \partial M, M],
 \qquad {\rm tr}\; K = {\rm tr}\; L =0.
 \end{equation}
 Explicitly, the matrix elements of $K$ and $L$ are of the form
 \begin{gather}
k_{11} = A^{-2} \left\{ ( \bar{w}_{1} \bar{\partial} w_{1}
 + \bar{w}_{2} \bar{\partial} w_{2})
 - ( w_{1} \bar{\partial} \bar{w}_{1} + w_{2} \bar{\partial} \bar{w}_{2})\right \},\nonumber\\
k_{12} = A^{-2} \left\{ |w_{1}|^2 \bar{\partial} w_{1} + w_{1} \bar{w}_{2}
 \bar{\partial} w_{2} - (\bar{\partial} w_{1} + w_{1}
 ( \bar{w}_{1} \bar{\partial} w_{1} + w_{1} \bar{\partial} \bar{w}_{1})\right.\nonumber\\
\left.\phantom{k_{12}=}{} + w_{2} ( \bar{w}_{2} \bar{\partial} w_{1} + w_{1} \bar{\partial}
 \bar{w}_{2})) \right\},
 \nonumber\\
k_{13} = A^{-2} \left\{ \bar{w}_{1} w_{2} \bar{\partial} w_{1}
 + |w_{2}|^2 \bar{\partial} w_{2} - (\bar{\partial} w_{2} + w_{1}
 ( \bar{w}_{1} \bar{\partial} w_{2} + w_{2} \bar{\partial} \bar{w}_{1})\right.\nonumber\\
\left.\phantom{k_{13}=}{}
 + w_{2} ( \bar{w}_{2} \bar{\partial} w_{2} + w_{2} \bar{\partial}
 \bar{w}_{2})) \right\},
 \nonumber\\
k_{21} = A^{-2}\left \{ \bar{\partial} \bar{w}_{1} + \bar{w}_{1}
 ( \bar{w}_{1} ( \bar{w}_{1} \bar{\partial} w_{1} + w_{1}
 \bar{\partial} \bar{w}_{1}) + \bar{w}_{2}
 ( \bar{w}_{1} \bar{\partial} w_{2} + w_{2}
  \bar{\partial} \bar{w}_{1})\right.
 \nonumber\\
\left.\phantom{k_{12}=}{}- (|w_{1}|^2 \bar{\partial} \bar{w}_{1}
 + \bar{w}_{1} w_{2} \bar{\partial} \bar{w}_{2}) \right\},
 \nonumber\\
k_{22} = A^{-2}\left\{ w_{1} \bar{\partial} \bar{w}_{1} \!+ w_{1} \bar{w}_{2}
 (\bar{w}_{1} \bar{\partial} w_{2}\! + w_{2} \bar{\partial} \bar{w}_{1})\!
 - (\bar{w}_{1} \bar{\partial} w_{1}\! + \bar{w}_{1} w_{2}
 (\bar{w}_{2} \bar{\partial} w_{1} \!+ w_{1} \bar{\partial} \bar{w}_{2}))\right \},
 \nonumber\\
k_{23} = A^{-2} \left\{ w_{2} \bar{\partial} \bar{w}_{1} + \bar{w}_{1} w_{2}
 ( \bar{w}_{1} \bar{\partial} w_{1} + w_{1} \bar{\partial} \bar{w}_{1} )
 + |w_{2}|^2 (\bar{w}_{1} \bar{\partial} w_{2} + w_{2} \bar{\partial}
 \bar{w}_{1})\right.
 \nonumber\\
\left.\phantom{k_{23}=}{}- \left[ \bar{w}_{1} \bar{\partial} w_{2} + |w_{1}|^2 (\bar{w}_{1}
 \bar{\partial} w_{2}
 + w_{2} \bar{\partial} \bar{w}_{1}) + \bar{w}_{1} w_{2} (\bar{w}_{2}
 \bar{\partial} w_{2} + w_{2} \bar{\partial} \bar{w}_{2}) \right] \right\},
 \nonumber\\
k_{31} = A^{-2} \left\{ \bar{\partial} \bar{w}_{2} + \bar{w}_{1}
 ( \bar{w}_{2} \bar{\partial} w_{1} + w_{1} \bar{\partial} \bar{w}_{2})
 + \bar{w}_{2} ( \bar{w}_{2} \bar{\partial} w_{2} + w_{2}
 \bar{\partial} \bar{w}_{2})\right.
 \nonumber\\
\left.\phantom{k_{31}=}{}- \left( w_{1} \bar{w}_{2} \bar{\partial} \bar{w}_{1} +
 | w_{2}|^2  \bar{\partial} \bar{w}_{2}\right) \right\},
 \nonumber\\
k_{32} = A^{-2} \left\{ w_{1} \bar{\partial} \bar{w}_{2} + |w_{1}|^2
 (\bar{w}_{2} \bar{\partial} w_{1} + w_{1} \bar{\partial} \bar{w}_{2})
 + w_{1} \bar{w}_{2} ( \bar{w}_{2} \bar{\partial} w_{2} + w_{2}
 \bar{\partial} \bar{w}_{2})\right.
 \nonumber\\
\left.\phantom{k_{32}=}{}- \left[ \bar{w}_{2} \bar{\partial} w_{1} + w_{1} \bar{w}_{2}
 ( \bar{w}_{1} \bar{\partial} w_{1} + w_{1} \bar{\partial} \bar{w}_{1} )
 + |w_{2}|^2 |\bar{w}_{2}|^2 (\bar{w}_{2} \bar{\partial} w_{1} + w_{1}
 \bar{\partial} \bar{w}_{2}) \right] \right\},
 \nonumber\\
k_{33} = A^{-2}\left\{ w_{2} \bar{\partial} \bar{w}_{2} + \bar{w}_{1} w_{2}
 ( \bar{w}_{2} \bar{\partial} w_{1} + w_{1} \bar{\partial} \bar{w}_{2})\right.
\nonumber\\
\left.\phantom{k_{33}=}{} - ( \bar{w}_{2} \bar{\partial} w_{2} + w_{1} \bar{w}_{2} ( \bar{w}_{1}
 \bar{\partial}w_{2} + w_{2}  \bar{\partial} \bar{w}_{1})) \right\},
 \label{4.42}
 \end{gather}
 and
 \begin{gather*}
l_{11} = A^{-2} \left\{ \bar{w}_{1} \partial w_{1} + \bar{w}_{2} \partial
 w_{2}
 - ( w_{1} \partial \bar{w}_{1} + w_{2} \partial \bar{w}_{2}) \right\},
 \nonumber\\
l_{12} = A^{-2} \left\{ |w_{1}|^2 \partial w_{1} + w_{1} \bar{w}_{2}
 \partial
 w_{2}
 - [ \partial w_{1} + w_{1} ( \bar{w}_{1} \partial w_{1}+ w_{1} \partial
 \bar{w}_{1})\right.
\nonumber\\
\left.\phantom{l_{12} =}{} + w_{2} (\bar{w}_{2} \partial w_{1} + w_{1} \partial \bar{w}_{2}) ] \right\},
 \nonumber\\
l_{13} = A^{-2} \left\{ \bar{w}_{1} w_{2} \partial w_{1} + |w_{2}|^2
 \partial w_{2} - [ \partial w_{2} + w_{1} (\bar{w}_{1} \partial w_{2} +
 w_{2} \partial
 \bar{w}_{1}) \right.\nonumber\\
\left.\phantom{l_{13}=}{}+ w_{2} ( \bar{w}_{2} \partial w_{2} + w_{2} \partial
 \bar{w}_{2})] \right\},
 \nonumber\\
l_{21} = A^{-2} \left\{ \partial \bar{w}_{1} + \bar{w}_{1} (\bar{w}_{1}
 \partial w_{1} +
 w_{1} \partial \bar{w}_{1}) + \bar{w}_{2} ( \bar{w}_{1} \partial w_{2} +
 w_{2}
 \partial \bar{w}_{1}) \right.\nonumber\\
\left.\phantom{l_{21}=}{}- ( |w_{1}|^2 \partial \bar{w}_{1} + \bar{w}_{1}
 w_{2}
 \partial \bar{w}_{2}) \right\},
 \nonumber\\
l_{22} = A^{-2} \left\{ w_{1} \partial \bar{w}_{1} + w_{1} \bar{w}_{2}
 ( \bar{w}_{1} \partial w_{2} + w_{2} \partial \bar{w}_{1})
 - [ \bar{w}_{1} \partial w_{1} + \bar{w}_{1} w_{2} (\bar{w}_{2}
 \partial w_{1} + w_{1} \partial \bar{w}_{2}) ] \right\},
 \nonumber\\
l_{23} = A^{-2} \left\{ w_{2} \partial \bar{w}_{1} + \bar{w}_{1} w_{2}
 ( \bar{w}_{1} \partial w_{1} + w_{1} \partial \bar{w}_{1}) + |w_{2}|^2
 (\bar{w}_{1} \partial w_{2} + w_{2} \partial \bar{w}_{1})\right.
 \nonumber\\
\left.\phantom{l_{23}=}{}- [ \bar{w}_{1} \partial w_{2} + |w_{1}|^2 (\bar{w}_{1} \partial w_{2}
 + w_{2} \partial \bar{w}_{1}) + \bar{w}_{1} w_{2} (\bar{w}_{2}
 \partial w_{2} + w_{2} \partial \bar{w}_{2})] \right\},
\end{gather*}
\begin{gather}
l_{31} = A^{-2} \left\{ \partial \bar{w}_{2} + \bar{w}_{1}
 ( \bar{w}_{2} \partial w_{1} + w_{1} \partial \bar{w}_{2})
 + \bar{w}_{2} ( \bar{w}_{2} \partial w_{2} + w_{2} \partial \bar{w}_{2})\right.
\nonumber\\
\left.\phantom{l_{31}=}{} - ( w_{1} \bar{w}_{2} \partial \bar{w}_{1} + |w_{2}|^2 \partial
 \bar{w}_{2} )  \right\},
 \nonumber\\
l_{32} = A^{-2} \left\{ w_{1} \partial \bar{w}_{2} + |w_{1}|^2
 ( \bar{w}_{2} \partial w_{1} + w_{1} \partial \bar{w}_{2})
 + w_{1} \bar{w}_{2} ( \bar{w}_{2} \partial w_{2} + w_{2} \partial
 \bar{w}_{2})\right.
 \nonumber\\
\left.\phantom{l_{32}=}{}- [ \bar{w}_{2} \partial w_{1} + w_{1} \bar{w}_{2} ( \bar{w}_{1}
 \partial w_{1} + w_{1} \partial \bar{w}_{1}) + |w_{2}|^2
 (\bar{w}_{2} \partial w_{1} + w_{1} \partial \bar{w}_{2})]  \right\},
 \nonumber\\
l_{33} = A^{-2} \left\{ w_{2} \partial \bar{w}_{2} + \bar{w}_{1} w_{2}
 ( \bar{w}_{2} \partial w_{1} + w_{1} \partial \bar{w}_{2}) \right.
\nonumber\\
\left.\phantom{l_{33}=}{}-
 (\bar{w}_{2} \partial w_{2} + w_{1} \bar{w}_{2} ( \bar{w}_{1}
 \partial w_{2} + w_{2} \partial \bar{w}_{1})) \right\},
 \label{4.43}
 \end{gather}
 respectively. Finally from equations (\ref{4.40}), (\ref{4.42})
 and (\ref{4.43}) we see that
 there exists only five independent conservation laws for $CP^2$
 sigma model (\ref{4.6})--(\ref{4.8}). Namely we have the following
 independent
 conserved quantities
 \begin{gather}
\partial \left\{ A^{-2} \left[ \bar{w}_{1} \bar{\partial} w_{1} + \bar{w}_{2}
 \bar{\partial} w_{2} - ( w_{1} \bar{\partial} \bar{w}_{1} + w_{2}
 \bar{\partial} w_{2} ) \right] \right\}\nonumber\\
\qquad {} +\bar{\partial} \left\{ A^{-2} \left[ \bar{w}_{1} \partial w_{1} + \bar{w}_{2}
 \partial w_{2}- ( w_{1} \partial \bar{w}_{1} + w_{2}
 \partial \bar{w}_{2}) \right] \right\} = 0, \nonumber\\
\partial \left\{ A^{-2} \left[ w_{1} \bar{\partial} \bar{w}_{1} + w_{1}
 \bar{w}_{2}
 ( \bar{w}_{1} \bar{\partial} w_{2} + w_{2} \bar{\partial} \bar{w}_{1} )
 - \bar{w}_{1} \bar{\partial} w_{1} - \bar{w}_{1} w_{2}
 ( \bar{w}_{2} \bar{\partial} w_{1} + w_{1} \bar{\partial} \bar{w}_{2}) \right]
 \right\}
 \nonumber\\
\qquad {}+\bar{\partial} \left\{ A^{-2} \left[ w_{1} \partial \bar{w}_{1} + w_{1}
 \bar{w}_{2}
 ( \bar{w}_{1} \partial w_{2} + w_{2} \partial \bar{w}_{1}) - \bar{w}_{1}
 \partial w_{1} \right.\right.
\nonumber\\
\qquad\left.\left.{}- \bar{w}_{1} w_{2} ( \bar{w}_{2} \partial w_{1} + w_{1}
 \partial \bar{w}_{2})\right] \right\}=0,\nonumber\\
\partial \left\{ A^{-2} \left[ |w_{1}|^2 \bar{\partial} w_{1} + w_{1} \bar{w}_{2}
 \bar{\partial} w_{2} - \bar{\partial} w_{1} - w_{1}
 ( \bar{w}_{1} \bar{\partial} w_{1} + w_{1} \bar{\partial} \bar{w}_{1})
\right.\right.
 \nonumber\\
\qquad \left.\left.{}
 - w_{2} (\bar{w}_{2} \bar{\partial} w_{1} + w_{1} \bar{\partial}
 \bar{w}_{2}) \right] \right\}
+ \bar{\partial}\left\{ A^{-2} \left[ |w_{1}|^2 \partial w_{1} + w_{1}
 \bar{w}_{2}
 \partial w_{2} - \partial w_{1}
\right.\right.
\nonumber\\
\qquad\left.\left.{} - w_{1} ( \bar{w}_{1} \partial w_{1} +
 w_{1}
 \partial \bar{w}_{1}) - w_{2} ( \bar{w}_{2} \partial w_{1} + w_{1}
 \partial \bar{w}_{2}) \right] \right\} \nonumber = 0,
 \nonumber\\
\partial \left\{ A^{-2}\left[\bar{w}_{1} w_{2} \bar{\partial} w_{1} -
 \bar{\partial} w_{2} - w_{1} ( \bar{w}_{1}
 \bar{\partial} w_{2} + w_{2} \bar{\partial} \bar{w}_{1}) - w_{2}^2
 \bar{\partial} \bar{w}_{2}\right]\right\}
 \nonumber\\
\qquad{}+ \bar{\partial}\left\{ A^{-2} \left[ \bar{w}_{1} w_{2} \partial w_{1} -
 \partial
 w_{2}
 - w_{1} ( \bar{w}_{1} \partial w_{2} + w_{2} \partial \bar{w}_{1}) -
 w_{2}^2
 \partial \bar{w}_{2} \right] \right\}  = 0,
 \nonumber\\
\partial \left\{ A^{-2} \left[ w_{2} \bar{\partial} \bar{w}_{1} + \bar{w}_{1}^2
 w_{2}
 \bar{\partial} w_{1} + |w_{2}|^2 w_{2} \bar{\partial} \bar{w}_{1}\right] \right.
\nonumber\\
\qquad\left.{}-  A^{-2}
\left[ \bar{w}_{1} \bar{\partial} w_{2} + |w_{1}|^2 \bar{w}_{1} \bar{\partial}
 w_{2}
 + \bar{w}_{1} w_{2}^2 \bar{\partial} \bar{w}_{2} \right]\right\}
 \nonumber\\
\qquad{}+ \bar{\partial} \left\{ A^{-2} \left[ w_{2} \partial \bar{w}_{1} + \bar{w}_{1}^2
 w_{2}
 \partial w_{1} + |w_{2}|^2 w_{2} \partial \bar{w}_{1} \right]\right.\nonumber\\
\qquad \left.{}
 - A^{-2} \left[ \bar{w}_{1} \partial w_{2} + |w_{1}|^2 \bar{w}_{1} \partial
 w_{2}
+\bar{w}_{1} w_{2}^2 \partial \bar{w}_{2} \right] \right\} = 0.
 \label{4.44}
\end{gather}

 Consequently, as a result of the conservation laws (\ref{4.44}) there
 exist eight
 real-valued functions $X^{i} (z, \bar{z})$, $i=1, \ldots, 8$ expressed in
 terms of
 functions $(w_{1},w_{2})$, i.e.\
\begin{gather*}
 X^1 =  \int_{C} A^{-2} \left\{ -\left[ \bar{w}_{1} \partial w_{1} +
 \bar{w}_{2}
 \partial w_{2} - (w_{1} \partial \bar{w}_{1} + w_{2} \partial
 \bar{w}_{2})\right] dz\right.\\
\left.\phantom{X^1=}{} + \left[ \bar{w}_{1} \bar{\partial} w_{1} + \bar{w}_{2} \bar{\partial} w_{2}
-(w_{1} \bar{\partial} \bar{w}_{1} + w_{2}
 \bar{\partial} \bar{w}_{2}) \right]  d \bar{z} \right\},\\
 X^2 =  \int_{C} A^{-2} \left\{ \left[ \left(1 + |w_{2}|^2 \right) (w_{1} \partial
 \bar{w}_{1} -\bar{w}_{1}
 \partial w_{1}) + |w_{1}|^2 (\bar{w}_{2} \partial w_{2} - w_{2} \partial
 \bar{w}_{2}) \right] dz\right.\\
\left.\phantom{X^2=}{}+ \left[ \left( 1 + |w_{2}|^2\right)(w_{1} \bar{\partial} \bar{w}_{1} - \bar{w}_{1}
 \bar{\partial} w_{1})
 + |w_{1}|^2 (\bar{w}_{2} \bar{\partial} w_{2} - w_{2} \bar{\partial}
 \bar{w}_{2})\right] d \bar{z} \right\},\\
 X^3 = i \int_{C} -A^{-2} \left\{ \left[ - \left(1 + \bar{w}_{1}^2 + |w_{2}|^2\right)
 \partial  w_{1}\right.\right.\\
\left.\phantom{X^3=}- \left(1 + w_{1}^2 + |w_{2}|^2\right) \partial \bar{w}_{1}
 + \bar{w}_{2} (w_{1} -  \bar{w}_{1})
 \partial w_{2} + w_{2} ( \bar{w}_{1} - w_{1}) \partial \bar{w}_{2} \right] dz
\\
\phantom{X^3=}{}+ \left[ - \left(1 + \bar{w}_{1}^2 + |w_{2}|^2 \right) \bar{\partial} w_{1} -
 \left( 1 + w_{1}^2 + |w_{2}|^2\right) \bar{\partial} \bar{w}_{1} + \bar{w}_{2}
 (w_{1} - \bar{w}_{1}) \bar{\partial} w_{2}  \right.\\
\left.\left.\phantom{X^3=}{}+ w_{2} (\bar{w}_{1} - w_{1} )
 \bar{\partial} \bar{w}_{2} \right] d \bar{z} \right\},
 \end{gather*}
 \begin{gather}
 X^4 = \int_{C} A^{-2} \left\{ \left[ \left(1 - \bar{w}_{1}^2 + |w_{2}|^2\right) \partial
 w_{1}\right.\right.
 \nonumber\\
\left.\phantom{X^4=}{}+ \left(-1 + w_{1}^2 - |w_{2}|^2\right) \partial \bar{w}_{1} - \bar{w}_{2} (w_{1}
 + \bar{w}_{1})
 \partial w_{2} + w_{2} ( w_{1} + \bar{w}_{1}) \partial \bar{w}_{2} \right] dz
 \nonumber\\
\phantom{X^4=}{}+ \left[ \left(-1 + \bar{w}_{1}^2 - |w_{2}|^2\right) \bar{\partial} w_{1}
 + \left(1 - w_{1}^2 + |w_{2}|^2\right) \bar{\partial} \bar{w}_{1} + \bar{w}_{2}
 (w_{1} + \bar{w}_{1}) \bar{\partial} w_{2} \right.\nonumber \\
\left.\left.\phantom{X^4=}{}- w_{2} (w_{1} + \bar{w}_{1})
 \bar{\partial} \bar{w}_{2} \right] d \bar{z}  \right\}, \nonumber\\
 X^5= i \int_{C} - A^{-2} \left\{ \left[ \bar{w}_{1} ( w_{2} - \bar{w}_{2})
 \partial w_{1}
 + w_{1} (\bar{w}_{2}- w_{2}) \partial \bar{w}_{1}\right.\right.
 \nonumber\\
\left.\phantom{X^5=}{}- \left(1 + |w_{1}|^2
 + \bar{w}_{2}^2\right) \partial w_{2} - \left(1 + |w_{1}|^2 + w_{2}^2\right) \partial
 \bar{w}_{2}\right] dz
 \nonumber\\
\phantom{X^5=}{}+ \left[ \bar{w}_{1} ( w_{2} - \bar{w}_{2}) \bar{\partial} w_{1} + w_{1}
 (\bar{w}_{2} - w_{2})
 \bar{\partial} \bar{w}_{1} - \left(1 + |w_{1}|^2 + \bar{w}_{2}^2\right)
 \bar{\partial} w_{2} \right.\nonumber \\
\left.\left.\phantom{X^5=}{}- \left(1 + |w_{1}|^2 + w_{2}^2\right)
\bar{\partial} \bar{w}_{2}\right] d \bar{z} \right\},
 \nonumber\\
 X^6 = \int_{C} A^{-2} \left\{ \left[ - \bar{w}_{1} (w_{2} + \bar{w}_{2}) \partial
 w_{1} + w_{1} (w_{2} + \bar{w}_{2}) \partial \bar{w}_{1}\right.\right.
 \nonumber\\
\left.\phantom{X^6=}{}+ \left(1 + |w_{1}|^2 - \bar{w}_{2}^2\right)
 \partial w_{2} - \left(1 + |w_{1}|^2 -w_{2}^2\right) \partial \bar{w}_{2} \right] dz
 \nonumber\\
\phantom{X^6=}{}+ \left[ \bar{w}_{1} (w_{2} + \bar{w}_{2}) \bar{\partial} w_{1} - w_{1}
 (w_{2} + \bar{w}_{2})
 \bar{\partial} \bar{w}_{2} - \left(1+ |w_{1}|^2 -\bar{w}_{2}^2\right) \bar{\partial}
 w_{2} \right.\nonumber \\
\left.\left.\phantom{X^6=}{}+\left(1+|w_1|^2-w_{2}^2\right) \bar{\partial} \bar{w}_{2}\right] d \bar{z}
\right\},
 \nonumber\\
 X^7 = i \int_{C} - A^{-2} \left\{ \left[ \bar{w}_{2} \left(1 + |w_{2}|^2\right) +
 \bar{w}_{1}^{2} w_{2} \right]
 \partial w_{1} + \left[ w_{2} (1 + |w_{2}|^2) + w_{1}^2 \bar{w}_{2} \right] \partial
 \bar{w}_{1}\right.
 \nonumber\\
\left.\phantom{X^7=}{}-\left[ \bar{w}_{1} \left( 1 + |w_{1}|^2 \right) + w_{1} \bar{w}_{2}^2\right]
\partial w_{2}
 - \left[ w_{1} \left(1 + |w_{1}|^2\right) + \bar{w}_{1} w_{2}^2 \right] \partial \bar{w}_{2} \right\}
 dz
 \nonumber\\
\phantom{X^7=}{}+ A^{-2} \left\{ \left[ \bar{w}_{2} \left(1 + |w_{2}|^2 \right) + \bar{w}_{1}^2 w_{2} \right]
 \bar{\partial} w_{1}
 + \left[ w_{2} \left(1 + |w_{2}|^2\right) + w_{1}^2 \bar{w}_{2}\right] \bar{\partial}
 \bar{w}_{1}\right.
 \nonumber\\
\left.\phantom{X^7=}{}- \left[ \bar{w}_{1}
\left(1 + |w_{1}|^2\right) + w_{1} \bar{w}_{2}^2 \right] \bar{\partial}
 w_{2}
 - \left[ w_{1} \left( 1 + |w_{1}|^2\right) + \bar{w}_{1} w_{2}^2\right] \bar{\partial}
 \bar{w}_{2} \right\} d \bar{z},
 \nonumber\\
 X^8 = \int_{C} A^{-2}\left\{\left[ \bar{w}_{2} \left(1 + |w_{2}|^2\right) - \bar{w}_{1}^2
 w_{2}\right] \partial w_{1}
 + \left[ - w_{2} \left(1 + |w_{2}|^2\right) + w_{1}^2 \bar{w}_{2}\right] \partial \bar{w}_{1}\right.
 \nonumber\\
\left.\phantom{X^8=}{}+
\left[ \bar{w}_{1}\left( 1 + |w_{1}|^2\right) - w_{1} \bar{w}_{2}^2\right] \partial w_{2} +
 \left[-w_{1}
\left(1 + |w_{1}|^2\right) + \bar{w}_{1} w_{2}^2 \right] \partial \bar{w}_{2} \right\} dz
 \nonumber\\
\phantom{X^8=}{}+ A^{-2} \left\{ \left[ - \bar{w}_{2} \left(1 + |w_{2}|^2\right) + \bar{w}_{1}^2 w_{2}\right]
 \bar{\partial} w_{1}
 +\left[ w_{2} \left(1 + |w_{2}|^2\right) - w_{1}^2 \bar{w}_{2}^2\right] \bar{\partial}
 \bar{w}_{1}\right. \nonumber\\
\left.\phantom{X^8=}{}+\left[ - \bar{w}_{1} \left( 1 + |w_{1}|^2\right) + w_{1} \bar{w}_{2}^2\right]
 \bar{\partial}  w_{2}
 + \left[ w_{1} \left(1 + |w_{1}|^2\right) - \bar{w}_{1} w_{2}^2 \right] \bar{\partial}
 \bar{w}_{2} \right\} d \bar{z}.
 \label{4.45}
 \end{gather}
 Note that by virtue of the conservation laws (\ref{4.44})
 for the $CP^2$ sigma model (\ref{4.6})--(\ref{4.8}) the r.h.s. in expression
 (\ref{4.45}) do not
 depend on the choice of the contour $C$ but only on its
 endpoints. This is due to the fact
 that (\ref{4.42}) are integrals of exact differentials
 of real valued functions. We
 identify
 the functions $X^{i} (z, \bar{z})$, $i=1, \ldots, 8$ with the coordinates
 of the radius
 vector
 \begin{equation}
 \boldsymbol{X} (z, \bar{z}) =
 \left( X^1 (z, \bar{z}), \ldots, X^8 (z, \bar{z}) \right),
 \end{equation}
 of a two-dimensional surface immersed into eight-dimensional Euclidean
 space
 ${\R}^8$. Substituting (\ref{4.10}) into (\ref{4.45}) we express the
 radius
 vector $\boldsymbol{X} (z, \bar{z})$ in terms of functions
 $(\varphi_{1}, \varphi_{2}, \psi_{1}, \psi_{2})$ and obtain
\begin{gather*}
 X^1 = 2 \int_{C} (\bar{\psi}_{1} \varphi_{1}+ \bar{\psi}_{2}
 \varphi_{2}) dz
 + ( \psi_{1} \bar{\varphi}_{1} + \psi_{2} \bar{\varphi}_{2}) d \bar{z},
\\
 X^2 = 2 \int_{C} \left\{ \frac{\bar{\psi}_{1}}{\varphi_{1}} \left[
 \varphi_{1}^2
 - \frac{\psi_{1}}{\bar{\varphi}_{1}} ( \bar{\psi}_{1} \varphi_{1}
 + \bar{\psi}_{2} \varphi_{2}) - \frac{\psi_{1}}{\bar{\varphi}_{1}}
 \partial \left(A^{-1}\right)\right]\right.
\end{gather*}
\begin{gather*}
\left.\phantom{X^2=}{}- \frac{A \bar{\varphi}_{2} \psi_{1}}{\Omega}
 \left[ \frac{J}{A^2} \frac{\psi_{2}}{\bar{\varphi}_{2}} + \left(\partial \left(A^{-1}\right)
 + \bar{\psi}_{1} \varphi_{1} + \bar{\psi}_{2} \varphi_{2}\right) \varphi_{2}^2\right]
 \right\} dz
 \nonumber\\
\phantom{X^2=}{}+ \left\{ \frac{\psi_{1}}{\bar{\varphi}_{1}} \left[ \bar{\varphi}_{1}^2
 - \frac{\bar{\psi}_{1}}{\varphi_{1}} ( \psi_{1} \bar{\varphi}_{1} +
 \psi_{2} \bar{\varphi}_{2}) - \frac{\bar{\psi}_{1}}{\varphi_{1}}
 \bar{\partial} \left(A^{-1}\right)\right]\right.
 \nonumber\\
\left.\phantom{X^2=}{}- \frac{A \varphi_{2} \bar{\psi}_{1}}
 {\bar{\Omega}} \left[ \frac{\bar{J}}{A^2}  \frac{\bar{\psi}_{2}}{\varphi_{2}}
 + \left( \bar{\partial}\left(A^{-1}\right) + \psi_{1} \bar{\varphi}_{1}
 + \psi_{2} \bar{\varphi}_{2}\right) \bar{\varphi}_{2}^2 \right] \right\}  d \bar{z},
 \nonumber\\
 X^3 = i \int_{C} - \left\{ - \frac{\bar{\psi}_{1}}{\varphi_{1}}
 \left[ \partial \left(A^{-1}\right) + 2 (\bar{\psi}_{1} \varphi_{1} + \bar{\psi}_{2}
 \varphi_{2})\right]\right.
 \nonumber\\
\phantom{X^3=}{}- \frac{A \bar{\varphi}_{1} \bar{\varphi}_{2}}{\Omega} \left[\frac{J}{A^2}
 \frac{\psi_{2}}{\bar{\varphi}_{2}} + \left( \partial\left(A^{-1}\right)
 + \bar{\psi}_{1} \varphi_{1} + \bar{\psi}_{2} \varphi_{2}\right) \varphi_{2}^2\right]
 \nonumber\\
\left.\phantom{X^3=}{}+ \left[ - \varphi_{1}^2 + \frac{\psi_{1}}{\bar{\varphi}_{1}}
 (\bar{\psi}_{1} \varphi_{1} + \bar{\psi}_{2} \varphi_{2}) +
 \frac{\psi_{1}}{\bar{\varphi}_{1}}
 \partial \left(A^{-1}\right)\right] \right\} dz
 \nonumber\\
\phantom{X^3=}{}+ \left\{ - \frac{\psi_{1}}{\bar{\varphi}_{1}} \left[ \bar{\partial} \left( A^{-1}\right) +
 2 (\psi_{1} \bar{\varphi}_{1} + \psi_{2} \bar{\varphi}_{2})\right]\right.
 \nonumber\\
\phantom{X^3=}{}- \frac{A \varphi_{1} \varphi_{2}}{\bar{\Omega}} \left[ \frac{\bar{J}}{A^2}
 \frac{\bar{\psi}_{2}}{\varphi_{2}} + \left( \bar{\partial} \left(A^{-1}\right) + \psi_{1}
 \bar{\varphi}_{1} + \psi_{2} \varphi_{2}\right) \bar{\varphi}_{2}^2 \right]
 \nonumber\\
\left.\phantom{X^3=}{}+ \left[ - \bar{\varphi}_{1}^2 + \frac{\bar{\psi}_{1}}{\varphi_{1}}
 (\psi_{1} \bar{\varphi}_{1} + \psi_{2} \bar{\varphi}_{2})
 + \frac{\bar{\psi}_{1}}{\varphi_{1}} \bar{\partial} \left(A^{-1}\right)\right] \right\}  d
 \bar{z},
 \nonumber\\
 X^4 = \int_{C} \left\{ - \frac{\bar{\psi}_{1}}{\varphi_{1}}
 \left[ \partial \left(A^{-1}\right) + 2 ( \bar{\psi}_{1} \varphi_{1} + \bar{\psi}_{2}
 \varphi_{2})\right]\right.
 \nonumber\\
\phantom{X^4=}{}- \frac{A \bar{\varphi}_{1} \bar{\varphi}_{2}}{\Omega} \left[ \frac{J}{A^2}
 \frac{\psi_{2}}{\bar{\varphi}_{2}} + \left( \partial\left(A^{-1}\right) + \bar{\psi}_{1}
 \varphi_{1} + \bar{\psi}_{2} \varphi_{2}\right) \varphi_{2}^2\right]
 \nonumber\\
\left.\phantom{X^4=}{}- \left[ - \varphi_{1}^2 + \frac{\psi_{1}}{\bar{\varphi}_{1}} (
 \bar{\psi}_{1}
 \varphi_{1} + \bar{\psi}_{2} \varphi_{2}) +
 \frac{\psi_{1}}{\bar{\varphi}_{1}}
 \partial \left(A^{-1}\right) \right] \right\} dz
 \nonumber\\
\phantom{X^4=}{}+\left\{ - \frac{\psi_{1}}{\bar{\varphi}_{1}} \left[ \bar{\partial} \left(A^{-1}\right) + 2
 (\psi_{1}
 \bar{\varphi}_{1} + \psi_{2} \bar{\varphi}_{2})\right]\right.
 \nonumber\\
\phantom{X^4=}{}- \frac{A \varphi_{1} \varphi_{2}}
 {\bar{\Omega}} \left[\frac{\bar{J}}{A^2} \frac{\bar{\psi}_{2}}{\varphi_{2}} +
\left( \bar{\partial} \left(A^{-1}\right) + \psi_{1} \bar{\varphi}_{1} + \psi_{2}
 \bar{\varphi}_{2}\right)
 \bar{\varphi}_{2}^2\right]
 \nonumber\\
\left.\phantom{X^4=}{}- \left[ - \bar{\varphi}_{1}^2 + \frac{\bar{\psi}_{1}}{\varphi_{1}}
 (\psi_{1} \bar{\varphi}_{1} + \psi_{2} \bar{\varphi}_{2})
 + \frac{\bar{\psi}_{1}}{\varphi_{1}} \bar{\partial} \left(A^{-1}\right)\right]\right\} d
 \bar{z},
 \nonumber\\
 X^5 = i \int_{C} - \left\{ - \frac{\bar{\psi}_{2}}{\varphi_{2}} \left[ \partial
 \left(A^{-1}\right)
 + 2 ( \bar{\psi}_{1} \varphi_{1} + \bar{\psi}_{2} \varphi_{2})\right]\right.
 \nonumber\\
\phantom{X^5=}{}+ \frac{A \bar{\varphi}_{1} \bar{\varphi}_{2}}{\Omega}
 \left[ \left( \partial \left(A^{-1}\right) + \bar{\psi}_{1} \varphi_{1} + \bar{\psi}_{2}
 \varphi_{2}\right)
 \varphi_{1}^2 + \frac{J}{A^2} \frac{\psi_{1}}{\bar{\varphi}_{1}} \right]
 \nonumber\\
\left.\phantom{X^5=}{}+ \left[ - \varphi_{2}^2 + \frac{\psi_{2}}{\bar{\varphi}_{2}}
 (\bar{\psi}_{1}
 \varphi_{1}
 + \bar{\psi}_{2} \varphi_{2}) + \frac{\psi_{2}}{\bar{\varphi}_{2}}
 \partial
 \left(A^{-1}\right) \right] \right\} dz
 \nonumber\\
\phantom{X^5=}{}+ \left\{ - \frac{\psi_{2}}{\bar{\varphi}_{2}} \left[ \bar{\partial} \left(A^{-1}\right)
 + 2 ( \psi_{1} \bar{\varphi}_{1} + \psi_{2} \bar{\varphi}_{2})\right]\right.
 \nonumber\\
\phantom{X^5=}{}+ \frac{A \varphi_{1} \varphi_{2}}{\bar{\Omega}} \left[ \left( \bar{\partial}
 \left(A^{-1}\right) + \psi_{1} \bar{\varphi}_{1} + \psi_{2}
 \bar{\varphi}_{2}\right) \bar{\varphi}_{1}^2 + \frac{\bar{J}}{A^2}
 \frac{\bar{\psi}_{1}}{\varphi_{1}}\right]\nonumber
\\
\left.\phantom{X^5=}{}+ \left[ - \bar{\varphi}_{2}^2 + \frac{\bar{\psi}_{2}}{\varphi_{2}}
 (\psi_{1} \bar{\varphi}_{1} + \psi_{2} \bar{\varphi}_{2})
 + \frac{\bar{\psi}_{2}}{\varphi_{2}} \bar{\partial} \left(A^{-1}\right) \right] \right \}  d  \bar{z},
\end{gather*}
\begin{gather*}
 X^6 = \int_{C} \left\{ - \frac{\bar{\psi}_{2}}{\varphi_{2}}
 \left[ \partial \left(A^{-1}\right) + 2 (\bar{\psi}_{1} \varphi_{1} + \bar{\psi}_{2}
 \varphi_{2}) \right]\right.
 \nonumber\\
\phantom{X^6=}{}+ \frac{A \bar{\varphi}_{1} \bar{\varphi}_{2}}{\Omega} \left[ \left( \partial \left(
 A^{-1}\right) +
 \bar{\psi}_{1} \varphi_{1} + \bar{\psi}_{2} \varphi_{2}\right) \varphi_{1}^2
 + \frac{J}{A^2} \frac{\psi_{1}}{\bar{\varphi}_{1}} \right]
 \nonumber\\
\left.\phantom{X^6=}{}- \left[ - \varphi_{2}^2 + \frac{\psi_{2}}{\bar{\varphi}_{2}}
 (\bar{\psi}_{1} \varphi_{1} + \bar{\psi}_{2} \varphi_{2}) +
 \frac{\psi_{2}}{\bar{\varphi}_{2}}
 \partial \left(A^{-1}\right) \right] \right\}  dz
 \nonumber\\
\phantom{X^6=}{}+ \left\{ - \frac{\psi_{2}}{\bar{\varphi}_{2}} \left[ \bar{\partial} \left( A^{-1}\right)
 + 2 ( \psi_{1} \bar{\varphi}_{1} + \psi_{2} \bar{\varphi}_{2})\right]\right.
 \nonumber\\
\phantom{X^6=}{}+ \frac{A \varphi_{1} \varphi_{2}}{\bar{\Omega}} \left[ \left( \bar{\partial}
 \left(A^{-1}\right)
 + \psi_{1} \bar{\varphi}_{1} + \psi_{2} \bar{\varphi}_{2}\right)
 \bar{\varphi}_{1}^2
 + \frac{\bar{J}}{A^2} \frac{\bar{\psi}_{1}}{\varphi_{1}}\right]
 \nonumber\\
\left.\phantom{X^6=}{}- \left[ - \bar{\varphi}_{2}^2 + \frac{\bar{\psi}_{2}}{\varphi_{2}} (
 \psi_{1}
 \bar{\varphi}_{1} + \psi_{2} \bar{\varphi}_{2}) +
 \frac{\bar{\psi}_{2}}{\varphi_{2}}
 \bar{\partial} \left(A^{-1}\right) \right] \right\}  d \bar{z},
 \nonumber\\
 X^7 = i \int_{C} - \left\{ \frac{\bar{\psi}_{2}}{\varphi_{2}} \left[
 \varphi_{1}^2
 - \frac{\psi_{1}}{\bar{\varphi}_{1}} ( \bar{\psi}_{1} \varphi_{1} +
 \bar{\psi}_{2} \varphi_{2})
 - \frac{\psi_{1}}{\bar{\varphi}_{1}} \partial \left(A^{-1}\right)\right]\right.
 \nonumber\\
\phantom{X^7=}{}+ \frac{A \bar{\varphi}_{2} \psi_{1}}{\Omega} \left[ \left(\partial \left(A^{-1}\right) +
 \bar{\psi}_{1}
 \varphi_{1} + \bar{\psi}_{2} \varphi_{2}\right) \varphi_{1}^2 + \frac{J}{A^2}
 \frac{\psi_{1}}{\bar{\varphi}_{1}}\right]
 \nonumber\\
\phantom{X^7=}{}+ \frac{A \bar{\varphi}_{1} \psi_{2}}{\Omega}\left[ \frac{J}{A^2}
 \frac{\psi_{2}}{\bar{\varphi}_{2}}
 +\left( \partial \left(A^{-1}\right) + \bar{\psi}_{1} \varphi_{1} + \bar{\psi}_{2}
 \varphi_{2}\right) \varphi_{2}^2\right]
 \nonumber\\
\left.\phantom{X^7=}{} + \frac{\bar{\psi}_{1}}{\varphi_{1}} \left[ - \varphi_{2}^2 +
 \frac{\psi_{2}}{\bar{\varphi}_{2}}
 ( \bar{\psi}_{1} \varphi_{1} + \bar{\psi}_{2} \varphi_{2}) +
 \frac{\psi_{2}}{\bar{\varphi}_{2}}
 \partial \left(A^{-1}\right)\right] \right\}  dz
 \nonumber\\
\phantom{X^7=}{}+ \left\{ \frac{\psi_{2}}{\bar{\varphi}_{2}} \left[ \bar{\varphi}_{1}^2 -
 \frac{\bar{\psi}_{1}}{\varphi_{1}}
 (\psi_{1} \bar{\varphi}_{1} + \psi_{2} \bar{\varphi}_{2}) -
 \frac{\bar{\psi}_{1}}{\varphi_{1}}
 \bar{\partial} \left(A^{-1}\right)\right]\right.
 \nonumber\\
\phantom{X^7=}{}+ \frac{A \varphi_{2} \bar{\psi}_{1}}{\Omega} \left[ \left(\bar{\partial}
 \left(A^{-1}\right)
 + \psi_{1} \bar{\varphi}_{1} + \psi_{2} \bar{\varphi}_{2}\right)
 \bar{\varphi}_{1}^2 + \frac{\bar{J}}{A^2}
 \frac{\bar{\psi}_{1}}{\varphi_{1}}\right]
 \nonumber\\
\phantom{X^7=}{}+ \frac{A \varphi_{1} \bar{\psi}_{2}}{\bar{\Omega}}\left[
 \frac{\bar{J}}{A^2}
 \frac{\bar{\psi}_{2}}
 {\varphi_{2}} + \left( \bar{\partial} \left(A^{-1}\right) + \psi_{1} \bar{\varphi}_{1} +
 \psi_{2}
 \bar{\varphi}_{2}\right) \bar{\varphi}_{2}^2\right]
 \nonumber\\
\left.\phantom{X^7=}{}+ \frac{\psi_{1}}{\bar{\varphi}_{1}}\left[ - \bar{\varphi}_{2}^2
 + \frac{\bar{\psi}_{2}}{\varphi_{2}} ( \psi_{1} \bar{\varphi}_{1} +
 \psi_{2} \bar{\varphi}_{2})
 + \frac{\bar{\psi}_{2}}{\varphi_{2}} \bar{\partial} \left(A^{-1}\right) \right] \right\} d
 \bar{z},
 \nonumber\\
 X^8 = \int_{C} \left\{ \frac{\bar{\psi}_{2}}{\varphi_{2}} \left[ \varphi_{1}^2 -
 \frac{\psi_{1}}{\bar{\varphi}_{1}}
 ( \bar{\psi}_{1} \varphi_{1} + \bar{\psi}_{2} \varphi_{2}) -
 \frac{\psi_{1}}{\bar{\varphi}_{1}}
 \partial \left(A^{-1}\right)\right ]\right.
 \nonumber\\
\phantom{X^8=}{}+ \frac{A \bar{\varphi}_{2} \psi_{1}}{\Omega} \left[ \left( \partial \left(A^{-1}\right)
 + \bar{\psi}_{1} \varphi_{1} + \bar{\psi}_{2} \varphi_{2}\right) \varphi_{1}^2
 + \frac{J}{A^2}
 \frac{\psi_{1}}{\bar{\varphi}_{1}}\right]
 \nonumber\\
\phantom{X^8=}{}- \frac{A \bar{\varphi}_{1} \psi_{2}}{\Omega} \left[ \frac{J}{A^2}
 \frac{\psi_{2}}{\bar{\varphi}_{2}}
 + \left( \partial \left(A^{-1}\right) + \bar{\psi}_{1} \varphi_{1} + \bar{\psi}_{2}
 \varphi_{2}\right)
 \varphi_{2}^2\right]
 \nonumber\\
\left.\phantom{X^8=}{}- \frac{\bar{\psi}_{1}}{\varphi_{1}} \left[ - \varphi_{2}^2 +
 \frac{\psi_{2}}{\bar{\varphi}_{2}} ( \bar{\psi}_{1} \varphi_{1} +
 \bar{\psi}_{2}
 \varphi_{2}) + \frac{\psi_{2}}{\bar{\varphi}_{2}} \partial \left(A^{-1}\right) \right] \right\}
 dz
 \nonumber\\
\phantom{X^8=}{}+ \left\{ \frac{\psi_{2}}{\bar{\varphi}_{2}} \left[ \bar{\varphi}_{1}^2 -
 \frac{\bar{\psi}_{1}}{\varphi_{1}}
 (\psi_{1} \bar{\varphi}_{1} + \psi_{2} \bar{\varphi}_{2} ) -
 \frac{\bar{\psi}_{1}}{\varphi_{1}}
 \bar{\partial} \left( A^{-1}\right)\right]\right.
 \nonumber\\
\phantom{X^8=}{}+ \frac{A \varphi_{2} \bar{\psi}_{1}}{\bar{\Omega}}
\left[ \left( \bar{\partial} \left(A^{-1}\right) + \psi_{1} \bar{\varphi}_{1} + \psi_{2}
 \bar{\varphi}_{2}\right)
 \bar{\varphi}_{1}^2 + \frac{\bar{J}}{A^2}
 \frac{\bar{\psi}_{1}}{\varphi_{1}}\right]\nonumber
\\
\phantom{X^8=}{}- \frac{A \varphi_{1} \bar{\psi}_{2}}{\bar{\Omega}} \left[
 \frac{\bar{J}}{A^2} \frac{\bar{\psi}_{2}}{\varphi_{2}}
 + \left( \bar{\partial} \left(A^{-1}\right) + \psi_{1} \bar{\varphi}_{1} + \psi_{2}
 \bar{\varphi}_{2}\right)
 \bar{\varphi}_{2}^2\right]
 \nonumber
\end{gather*}
\begin{gather}
\left.\phantom{X^8=}{}- \frac{\psi_{1}}{\bar{\varphi}_{1}} \left[ - \bar{\varphi}_{2}^2
 + \frac{\bar{\psi}_{2}}{\varphi_{2}} (\psi_{1} \bar{\varphi}_{1}
 + \psi_{2} \bar{\varphi}_{2})
 + \frac{\bar{\psi}_{2}}{\varphi_{2}} \bar{\partial} \left(A^{-1}\right) \right] \right\} d
 \bar{z},
 \label{4.47}
 \end{gather}
 where we have introduced the following notation
 \begin{equation}\label{4.48}
 \Omega = \varphi_{1} \psi_{2} | \varphi_{1}|^2 - \varphi_{2} \psi_{1} |
 \varphi_{2}|^2,
 \qquad
 \bar{\Omega} = \bar{\varphi}_{1} \bar{\psi}_{2} | \varphi_{1}|^2
 - \bar{\varphi}_{2} \bar{\psi}_{1} | \varphi_{2}|^2.
 \end{equation}

 Note that when  $J$
 vanishes, as can be checked, the position vector $X$ given by
 (\ref{4.47}) obeys the following relations
 \begin{equation}\label{4.49}
 ( \partial X, \partial X) = 0,
 \end{equation}
 and
 \begin{equation}\label{4.50}
 ( \partial \bar{\partial} X, \partial \bar{\partial} X) = ( \partial X,
 \bar{\partial} X)^{2} \neq 0,
 \end{equation}
 whenever the equations of the $CP^2$ sigma model (\ref{4.6})--(\ref{4.8})
 are satisfied. The
 explicit form of (\ref{4.50}), when written in terms of $w$'s or $\psi$'s
 and
 $\phi$'s, is very complicated and so we shall not reproduce it
 here. As a consequence of (\ref{4.49}) and (\ref{4.50}) the components of
 induced
 metric are
 \begin{equation}\label{4.51}
 g_{zz} = g_{\bar{z}\bar{z}} = 0,\qquad g_{z\bar{z}} \neq 0
 \end{equation}
 and the norm of the mean curvature vector $\bar{H} = (g_{z\bar{z}})^{-1}
 \partial
 \bar{\partial} X$, as expected, is equal to one, i.e.\ $|\bar{H}|^{2} =
 1$.
 We would like to note that when $w_{i}$ are  holomorphic functions then
 $J
 = 0$
 and formulae (\ref{4.45}) define a surface on $SU(3)$. Then using
 expressions in~\cite{Fok} we can calculate, in a closed form, all
 geometric characteristics of a given
 surface.

 Thus we have proved that
    the conformal immersions of CMC-surfaces into ${\R}^8$ are
 determined by formulae (\ref{4.45}) or (\ref{4.47}), where the complex
 functions
 $w_i$  obey the equations of $CP^2$ sigma model
 (\ref{4.6})--(\ref{4.8}), (or complex
 functions $\psi_i$ and $\phi_i$ obey the first order system
 (\ref{4.1})--(\ref{4.4})) and  $J$ given by (\ref{4.29}) (or (\ref{4.31}))
 vanishes.

 \section{Examples and applications}

In this section, based on our results of previous
 sections  we  construct certain classes of
 two-dimensional CMC-surfaces immersed into ${\R}^8$. For this purpose we use
 the
 $CP^2$ sigma model defined over $S^2$. Note that for such a model all
 solutions of
 the Euler--Lagrange equations (\ref{4.6})--(\ref{4.8}) are well known~\cite{book11}.
 Under the
 requirement of the finiteness of the action they split into three
 separate
 classes, i.e.\ analytic (i.e.\ $w_i = w_i(z)$), antianalytic (i.e.\ $w_i =
 w_i(\bar{z})$) and mixed ones. The latter ones can be determined from
 either
 the
 holomorphic or antiholomorphic functions by the following procedure.

   Consider three arbitrary holomorphic functions $f_i = f_i(z)$ and
 define
 for each pair the Wronskian
  \begin{equation}\label{5.1}
 F_{ij} = f_i\partial f_j -
 f_j \partial f_i,
 \qquad
 \bar {\partial}f_i = 0,
 \qquad
 i,j = 1,2,3
 \end{equation}
 Next determine three complex valued functions
 \begin{equation}\label{5.2}
 g_i = \sum^3_{k\neq i} \bar{f_k}
 (f_k\partial f_i -
 f_i \partial f_k),\qquad
 i = 1,2,3
 \end{equation}
 Then the mixed solutions $w_i$ of $CP^2$ sigma model
 (\ref{4.6})--(\ref{4.8}) can be
 determined
 as ratios of the components of $g_i$, i.e.\
 \begin{equation}\label{5.3}
 w_1 = \frac{g_1}{g_3},
 \qquad
 w_2 = \frac{g_2}{g_3},
 \qquad
 g_3\neq 0.
 \end{equation}
 Alternatively, similar class of solutions can be obtained when we
 consider
 three arbitrary antiholomorphic functions $f_i = f_i(\bar{z})$ and
 construct $g_i$ in the same way as above, but using~$\bar{\partial}$
 instead of $\partial$ in the equations~(\ref{5.2}).

  Now, let us discuss some classes of CMC-surfaces in ${\R}^n$ which can be
 obtained directly by applying the Weierstrass representation (\ref{4.45})
 and~(\ref{5.3}).

 \medskip

   {\bf 1}. One of the simplest solutions which corresponds to
 the analytic
 choice of functions is
 \begin{equation}\label{5.4}
 w_1 = z,
 \qquad
 w_2 = 1,
 \qquad
 A = 2 + |z|^2,
 \qquad
 J = 0.
 \end{equation}

 Using the $CP^2$ representation (\ref{4.45}) we can find that the
 associated CMC-surface is immersed in ${\R}^3$ and is given in a polar
 coordinates $\left(r = \left(x^2+y^2\right)^{1/2}, \varphi\right)$ by
 \begin{gather}
X^1 = -2 X^2 = 2 X^6 = 2 \sqrt2 \left(2+r^2\right)^{-1},
 \qquad
 X^3 = - X^7 = -2 r \left(2+r^2\right)^{-1} \sin\varphi,
 \nonumber\\
X^5 = 0, \qquad
  X^4 = X^8 = 2 r \left(2+r^2\right)^{-1} \cos\varphi.
 \label{5.5}
 \end{gather}
 The metric is conformally flat
 \begin{equation}\label{5.6}
 I = \frac{2}{\left(2+r^2\right)^2}\,\left(dr^2 + r^2 d\varphi^2\right)
 \end{equation}
 This case corresponds to the immersion of the $CP^2$ model into
 the $CP^1$ model.

 \medskip

   {\bf 2}. Another class of two-soliton solutions of the $CP^2$ model
 (\ref{4.6})--(\ref{4.8})
 is determined, for example, by two analytic functions; i.e.\ we can take,
 for example:
 \begin{equation}\label{5.7}
 w_1 = z^2,
 \qquad
 w_2 = \sqrt2 z,
 \qquad
 A = \left(1+|z|^2\right)^2,
 \qquad
 J = 0.
 \end{equation}

 Integrating formulae (\ref{4.45}) we obtain the associated CMC-surface
 which can be written in a polar coordinates as follows
 \begin{gather}
X^1 = 2 \left(1+r^2\right)^{-2},
 \qquad
 X^2 = -2 \left(1+2r^2\right)\left(1+r^2\right)^{-2},
 \nonumber\\
 X^3 = -2 r^2 \left(1+r^2\right)^{-2} \sin2\varphi,
 \qquad
 X^4 = 2 r^2 \left(1+r^2\right)^{-2} \cos2\varphi,
 \nonumber\\
 X^5= -2\sqrt2 r \left(1+r^2\right)^{-2} \sin\varphi,
 \qquad
X^6 = 2\sqrt2 r \left(1+r^2\right)^{-2} \cos\varphi,
 \nonumber\\
 X^7 = \sqrt2 \left(2r^2 -3\right)r\left(1+r^2\right)^{-2}\sin\varphi,
 \nonumber\\
X^8 = \sqrt2 \left(2r^2 -3\right)r\left(1+r^2\right)^{-2}\cos\varphi,
 \label{5.8}
 \end{gather}
 The corresponding first fundamental form is conformal
 \begin{equation}\label{5.9}
 I = \frac{2}{\left(1+r^2\right)^2}\,\left(dr^2 + r^2 d\varphi^2\right).
 \end{equation}

 \medskip

   {\bf 3}. A class of mixed solutions of the $CP^2$ model
 (\ref{4.2}) is represented by
 \begin{gather}
w_1 = \frac{\bar{z} + z}{1 - |z|^2},
 \qquad
 w_2 = \frac{\bar{z} - z}{1 - |z|^2},
 \qquad
A = \left(\frac{1 + |z|^2}{1 - |z|^2}\right)^2,
 \qquad
 J = 0.
 \label{5.10}
 \end{gather}
 From (\ref{4.32}) and using (\ref{5.2}) we obtain the expression for the
 associated surface which can be written in polar coordinates as follows
 \begin{gather}
X^1 = X^2 = X^4 = X^5 = X^7 = 0,
  \nonumber\\
 X^3 = \frac{-4 r}{1 + r^2}\, \sin\varphi,
 \qquad
 X^6 = \frac{-4 r}{1 +r^2}\, \cos\varphi,
 \qquad
 X^8 = \frac{4}{1 +r^2}\,
 \label{5.11}
 \end{gather}
 Hence this CMC-surface is really immersed in ${\R}^3$. The metric is
 conformal
 \begin{equation}\label{5.12}
 I = \frac{16}{\left(1 +r^2\right)^2}\,\left(dr^2 + r^2 d\varphi^2\right).
 \end{equation}
 This case corresponds to the immersion of the $CP^2$ model into the
 $CP^1$ model.

  If we take a more complicated example of this class,
 say,
 \begin{gather}
w_1 = \frac{-\bar{z}\left(3+2|z|^2\right)}{3 - |z|^4},
 \qquad
 w_2 = \frac{\sqrt{3}z\left(2+|z|^2\right)}{3 - |z|^4},
  \nonumber\\
 A = \frac{\left(1 + |z|^2\right)\left(|z|^6+6|z|^4+12|z|^2+9\right)}
{\left(3 - |z|^4\right)^2},\qquad
  J = 0.
 \label{5.113}
 \end{gather}
 then the expressions become very complicated but, in this case, we do
 have a genuine $CP^2$ solution.

 \medskip

 {\bf 4}. An interesting class of meron-like solutions of the $CP^1$
 model~(\ref{3.5}) is given by
 \begin{equation}\label{5.13}
 w = \left(\frac{z}{\bar{z}}\right)^{\beta},
 \qquad
 A = 2,
 \qquad
 J = \frac{-\beta^2}{2 z^2}.
 \end{equation}
 Here $\beta$ can be an integer of a half-integer, with $2\beta$ being
 the meron number. Note also that all merons are located at $z=0$
 and so this solution is defined on ${\R}^2\backslash\{0\}$.

 Then using the transformation
 (\ref{3.7}) we find that the solution of the GW system (\ref{3.1}) is
 given
 by
 \begin{equation}\label{5.14}
 \psi_1 = \frac{\sqrt{\beta}}{2\sqrt{\bar z}}\left(\frac{z}{\bar z}\right)
 \sp{\frac{1}{2}},
 \qquad
 \psi_2 = \frac{\sqrt{\beta}}{2\sqrt{ z}}\left(\frac{z}{\bar z}\right)
 \sp{\frac{1}{2}},
 \qquad
 p = \frac{\beta}{2|z|}.
 \end{equation}
 The associated surface  is obtained by integrating
 formulae~(\ref{3.20}) and we find
 \begin{equation}\label{5.15}
 X_1 = \frac{-3}{2} \sin 2 \beta\varphi,
 \qquad
 X_2 = \frac{-3}{2} \cos 2 \beta\varphi,
 \qquad
 X_3 = \frac{\beta}{2} \ln r.
 \end{equation}
 Thus the CMC-surface is a cylinder which is covered $2\beta$ times.
  Of course, the Gauss and mean curvatures are
 \begin{equation} \label{5.16}
 K = 0,
 \qquad
 H = 1.
 \end{equation}

 We note that our procedure of using the `multimeron' solution
 of the $CP^1$ model to construct our CMC-surface has effectively
 involved mapping  ${\R}^2\backslash\{0\}$ onto the equator of the
 unit sphere (merons) and then turning this circle into an infinite
 cylinder
 (based on this circle). This was done by effectively `undoing' the radial
  projection of the previous map.
 Hence the two singular points of the multi-meron
 configuration (\ref{5.13}) i.e.\ ($z=0$ and the point at $\infty$)
 have got mapped at the `ends' of the cylinder.

 \section{Final remarks and future developments}

   In this paper we have demonstrated links between the $CP^1$ and $CP^2$
 sigma models and
 Weierstrass representations for two-dimensional surfaces immersed into
 Euclidean spaces~${\R}^3$ and ${\R}^8$,
respectively. These links enabled us to
 present an algorithm for the construction of CMC-surfaces immersed into
 ${\R}^n$. This new approach has been tested in Section~5. It has proved to
 be
 effective, as we were able to reproduce easily the known results which,
 before,
 were
 obtained  by much  more complicated procedures. Its potential for
 providing new meaningful results has been exhibited in the case of mixed
 solutions of the $CP^2$ sigma model leading to new interesting surfaces
 in
 ${\R}^8$.

   The analytic method of the construction of the CMC-surfaces
 immersed into ${\R}^n$
 presented here is limited by several assumptions. For example we
 have studied
 only  $CP^1$ and $CP^2$ models defined over $S^2$ (i.e.\ to consider
 $J\ne 0$)
 The question the arises as to whether our approach can be extended to
 higher
   $CP^N$  models and to Weierstrass systems descri\-bing
 surfaces
 immersed in multi-dimensional Euclidean and pseudo-Riemannian spaces.
 If this is the case
 our approach may provide new classes of solutions and consequently new
 classes
 of surfaces in these multi-dimensional spaces. Other requirement of the
 proposed method, worth investigating further is
 the $CP^N$ models involving maps
 from ${\R}^2$ (not necessarly $S^2$). We can expect that
 taking $J\ne 0$  can broaden the applicability of
 our approach.

 Finally, it is worth noting that the CMC-surfaces can be used ``in
 reverse'' to address certain physical problems. Namely, we sometimes know
 the analytical description of CMC-surfaces in a physical system for
 which analytic models are not fully developed. Using our approach we can,
 perhaps,
 select an appropriate sigma model corresponding to the given Weierstrass
 representation and characterise the class of equations describing the
 physical phenomena in question. This was attempted successfully for
 Weierstrass representation for CMC-surfaces in 3-dimensional Euclidean
 space~\cite{book13}, but not to our knowledge for multi-dimensional
 spaces. These and other questions will be addressed in future work.

 \subsection*{Acknowledgements}

 AMG thanks A~Strasburger (University of Warsaw) for
 helpful and interesting discussions on the topic of this paper. WJZ would
 like to thank the CRM, Universite de Montreal and Center for
 Theoretical physics, MIT for the support of his stay and their
 hospitality. AMG would like to thank the University of Durham for the
 award
 of Alan Richards fellowship to him that allowed him to spend two terms in
 Durham during the academic years 2001 and 2002. Partial support for AMG's
 work was provided also by a grant
 from NSERC of Canada and the Fonds FCAR du Quebec.

 \label{grundland-lastpage}


\begin{thebibliography}{88}
\small

 \bibitem{Enneper}
 Enneper A,
{\it  Nachr. Konigl. Gesell. Wissensch. Georg - Augustus Univ. Gottingen}
 {\bf 12} (1868), 258--277.

 \bibitem{Weierstrass}
Weierstrass K, Fortsetzung der Untersuchung uber die
 Minimalflachen, Mathematische Werke,
 Vol.~3, Veragsbuch-handlung, Hillesheim, 1866, 219--248.

 \bibitem{book1}
Oserman R, A Survey of Minimal Surfaces,
 Dover, New York, 1996.

 \bibitem{K}
Konopelchenko B,
 Induced Surfaces and Their Integrable Dynamics, {\it Stud. Appl. Maths.} {\bf 96}
(1996), 9--51.

 \bibitem{book2}
Bobenko A~I,
 Surfaces in Terms of 2 by 2 Matrices.Old and New Integrable Cases,
 in Harmonic Maps and Integrable System, Editors: Fordy~A~P and Wood~J~C,
 Vieweg, Wiesbaden, 1994, 193--202.

 \bibitem{book3}
Gross~D~G, Pope C~N and Weinberg~S (Editors),
 Two-Dimensional Quantum Gravity and Random Surfaces,
  World Scientific, Singapore, 1992.

 \bibitem{Carroll}
Carroll R and Konopelchenko B,
 Generalised Weierstrass--Enneper Inducing Conformal Immersions and
 Gravity,
 {\it Int. J. Mod. Phys.} {\bf A11}, Nr.~7 (1996), 1183--1216.

 \bibitem{Land}
Konopelchenko B and Landolfi~G,
 On Classical String Configurations,
{\it Modern Phys. Letts.} {\bf 12} (1997), 3161--3179.

 \bibitem{Vis}
Viswanathan K and Parthasarathy R,
 {\it  Phys. Rev. D} {\bf D51} (1995), 5830--5879.

 \bibitem{book5}
Nelson D, Piran T and Weinberg S,
 Statistical Mechanics of Membranes and Surfaces,
 World Scientific, Singapore, 1989.

 \bibitem{book6}
Amit D,
 Field Theory, the Renormalization Group and Critical Phenomena,
 McGraw-Hill, New York, 1978.

 \bibitem{book7}
Rozdestvenski B~I and Yanenko N~N,
 Systems of Quasilinear Equations and their Applications to Gas Dynamics,
 AMS. Transl., Vol.~55, Providence, RI, 1983.

 \bibitem{Z}
Ou-Yang Z~O, Liu J and Xie Y~Z,
 Geometric Methods in the Elastic Theory of Membranes in Liquid Crystal
 Phases,  World Scientific, Singapore, 1999.

 \bibitem{Can}
Canham P~B, {\it  J.Theor. Biol.} {\bf 26} (1970), 61--81.

 \bibitem{Hel}
Helfrich W, {\it  Z. Naturforsch.} {\bf 28} (1973), 693--709.

 \bibitem{book8}
 Safran~S~A,
 Statistical Thermodynamics of Surfaces. Interfaces and Membranes,
 Addison Wesley, Reading, 1994.

 \bibitem{book9}
Eisenhart L~P,
 Treatise on the Differential Geometry of Curves and Surfaces,
 Dover, New York, 1909.

 \bibitem{book12}
Willmore T~J,
 Total Curvature in Riemannian Geometry,
 Ellis Horwood, New York, 1982.

 \bibitem{Tai}
Konopelchenko B and Taimanov~I,
 Constant Mean Curvature Surfaces via an Integrable Dynamical System,
{\it  J. Phys.} {\bf A29} (1996), 1261--1265.

 \bibitem{Land2}
Konopelchenko B and Landolfi G,
 Generalised Weierstrass Representation for Surfaces in Multi-Dimensional
 Riemann Spaces,
{\it J. Geom. Phys.} {\bf 29} (1999), 314--333.

 \bibitem{Bra3}
Bracken P and Grundland A~M,
 Symmetry Properties and Explicit Solutions of the Gene\-ra\-li\-sed Weierstrass
 System, {\it J. Math. Phys.} {\bf 42}, Nr.~3 (2001), 1250--1282.

 \bibitem{Ken}
Kenmotsu K,
 Weierstrass Formula for Surfaces of Prescribed Mean Curvature,
 {\it Math. Ann.} {\bf 245} (1979), 89--99.

 \bibitem{book10}
Zakharov V~E and Mikhailov A~V,
 Relativistically invariant two dimensional models in field theory which are integrable by means of the inverse scatering problem method,
{\it Sov. Phys.-JETP.} {\bf 47} (1979), 1017--1059.

 \bibitem{book11}
Zakrzewski W~J,
 Low Dimensional Sigma Models,
  Adam Hilger, Bristol, 1989.

 \bibitem{Abe}
Abe K and Erbacher J,
 Isometric Immersions with the Same Gauss Map,
{\it  Math. Ann.} {\bf 215} (1975), 197--201.

 \bibitem{Hof}
Hoffman~D and Osserman R,
 The Gauss Map of Surfaces in ${\R}^n$,
{\it  J. Diff. Geom.} {\bf 18} (1983), 733--754.

 \bibitem{Fer}
Ferapontov~E~V and Grundland~A~M,
 Links Between Different Analytic Descriptions of Constant Mean Curvature
 Surfaces, {\it  J. Nonlinear Math. Phys.} {\bf 7} (2000),  14--21.

 \bibitem{G}
Grundland A~M and Zakrzewski W~J,
 The Weierstrass Representation for Surfaces Immersed into ${\R}^8$ and
 $CP^2$
 Maps,
 {\it J. Math. Phys.} {\bf 43} (2002), 3352--3362.

 \bibitem{Fok}
Fokas A~S and Gelfand~I~M,
 Surfaces on Lie Groups, on Lie Algebras and Their Integrability,
{\it  Comm. Math. Phys.} {\bf 177} (1996), 203--220.

 \bibitem{book13}
Grosse-Brauckman K and Polthier K,
 Constant Mean Curvature Surfaces Derived from Delauney's and Wente's
 Examples, in Visualization and Mathematics: Experiments, Simulations and
 Environments, Editors: Hege~H~C and Polthier~K, Springer-Verlag, Berlin,
 1997, 386--419.
 \end{thebibliography}
 \end{document}